\documentclass[12pt]{article}
\usepackage{amssymb}
\usepackage{amsmath}
\usepackage{array}
\usepackage[russian]{babel}
\usepackage[pdftex]{graphicx}
\usepackage[X2,T2A]{fontenc}
\usepackage[cp1251]{inputenc}
\usepackage{mathrsfs}
\usepackage{longtable}
\usepackage{mathrsfs}
\usepackage{textcase}

\DeclareMathOperator{\RRe}{Re} \DeclareMathOperator{\IIm}{Im}

\multlinegap=0em \DeclareFontFamily{T1}{msb}{}
\DeclareFontShape{T1}{msb}{m}{ol}{<5> <6> <7> <8> <9> gen * msbm
<10> <10.95> <12> <14.4> <17.28> <20.74> <24.88> msbm10}{}
\DeclareSymbolFont{AMSb}{T1}{msb}{m}{ol} \multlinegap=0em

\renewcommand{\S}{\mathhexbox278}

\begin{document}

\begin{center}
{\rmfamily\bfseries\normalsize An elementary proof of some Ramanujan-type identities}
\end{center}

\begin{center}
{\normalsize M.A.~Korolev\footnote{Steklov Mathematical Institute of Russian Academy of Sciences. E-mail: \texttt{korolevma@mi-ras.ru}, \texttt{hardy\_ramanujan@mail.ru}}}
\end{center}

\vspace{1cm}

\fontsize{11}{12pt}\selectfont

The well-known Euler formula expresses the value of the Riemann zeta function $\zeta(s)$ at even integers as a product of the corresponding even power of $\pi$ and a rational number. In the first half of the 19th century, Cauchy found that the value $\zeta(3)$ could be expressed as the sum of some rational multiple of $\pi^{3}$ and a rapidly convergent infinite series. Later, Cauchy's formula was rediscovered and generalized by Lerch, Ramanujan, and others.
In this paper, we present an elementary proof of certain identities that express the squares of the Riemann zeta function at integer points as linear combinations of zeta values at integer points and some series involving hyperbolic functions and the digamma function. Cauchy's formula appears to be a special case of our relations.
We also prove a general identity involving double series whose coefficients depend on an arbitrary arithmetical function and its Dirichlet convolution. Finally, we study the properties of the special coefficients that appear naturally in this context.

Bibliography: 9 titles.

\vspace{1cm}

\fontsize{12}{15pt}\selectfont

\textbf{\S 1. Introduction.}
\vspace{0.5cm}

The main goal of this note is to apply a general but elementary identity, established in \cite{Daniyarkhodzhaev_Korolev_2021}, to derive some relations that connect the squares of the Riemann zeta-function $\zeta(z)$ at integer values with the digamma function, Bernoulli numbers and other related objects.

The following lemma is a special case of the general assertion from \cite{Daniyarkhodzhaev_Korolev_2021}:
\vspace{0.3cm}

\noindent
\textsc{Lemma 1.} \emph{Suppose that $\{a_{n}\}$, $\{b_{n}\}$, $n = 1,2,\ldots, N$, are arbitrary sequences such that $b_{n}>0$ for any $n$. Then the following identity holds:}
\begin{equation}\label{lab_01}
\sum\limits_{m,n=1}^{N}a_{m}a_{n}\,\frac{b_{m}}{b_{m}+b_{n}} = \frac{1}{2}\biggl(\,\sum\limits_{n=1}^{N}a_{n}\biggr)^{2}.
\end{equation}
For the reader's convenience, we give here a short proof of (\ref{lab_01}):
\begin{multline*}
\biggl(\,\sum\limits_{n=1}^{N}a_{n}\biggr)^{2} = \sum\limits_{m,n=1}^{N}a_{m}a_{n} = \sum\limits_{m,n=1}^{N}a_{m}a_{n}\frac{b_{m}+b_{n}}{b_{m}+b_{n}} = \\
= \sum\limits_{m,n=1}^{N}a_{m}a_{n}\,\frac{b_{m}}{b_{m}+b_{n}} + \sum\limits_{m,n=1}^{N}a_{m}a_{n}\,\frac{b_{n}}{b_{m}+b_{n}}  = 2\sum\limits_{m,n=1}^{N}a_{m}a_{n}\,\frac{b_{m}}{b_{m}+b_{n}}.\quad \square
\end{multline*}
In \cite{Daniyarkhodzhaev_Korolev_2021}, relations of this type were used to establish several curious identities for the values of the Riemann zeta function. For example, it was shown that
\[
\sum\limits_{m,n=1}^{+\infty}\frac{1}{m^{2}(m^{2}+n^{2})} = \frac{1}{2}\,\zeta^{2}(2) = \frac{\pi^{4}}{72},\quad
\sum\limits_{m,n=1}^{+\infty}\frac{1}{m^{3}(m^{3}+n^{3})} = \frac{1}{2}\,\zeta^{2}(3).
\]
In what follows, we provide some applications of (\ref{lab_01}) that seem to be new.
\vspace{0.5cm}

\noindent
\textsc{Notations.} The symbols $k,\ell,m,n,q,r,s$ denote integers. For $k\geqslant 1$ and $r\geqslant 0$ we set
\[
\varepsilon_{r} = \varepsilon_{r}(k) = \exp{\biggl(\frac{2\pi i(r+1/2)}{2k}\biggr)},\quad
\omega_{r} = \omega_{r}(k) = \exp{\biggl(\frac{2\pi i (r+1/2)}{2k+1}\biggr)}.
\]
Next, for positive real $w$, we define
\begin{align*}
& a_{k}(w) = \frac{1}{k}\sum\limits_{r=0}^{k-1}\pi\varepsilon_{r}\ctg{(\pi\varepsilon_{r}w)},\quad
b_{k,\ell}(w) = \frac{(-1)^{\ell+1}}{k}\sum\limits_{r=0}^{2k-1}\varepsilon_{r}^{\ell+1}\psi(\varepsilon_{r}w),\quad
b_{k}(w) = b_{k,1}(w),\\
& c_{k}(w) = \frac{1}{2k+1}\sum\limits_{r=0}^{2k}\omega_{r}\psi(-\omega_{r}w),
\end{align*}
where $\psi(z) =\Gamma'(z)/\Gamma(z)$ is the digamma-function. We also use standard notation for the basic arithmetic functions: the M\"{o}bius function $\mu(n)$, the divisor function $\tau(n)$ and its multidimensional analogue $\tau_{\nu}(n)$, the Liouville function $\lambda(n)$, the von Mangoldt function $\Lambda(n)$. The symbol $\sigma_{k}(n)$ denotes the sum of the $k$-th powers of the divisors of $n$, and $\omega(n)$ denotes the number of distinct prime divisors of $n$.
\vspace{0.5cm}

\textbf{\S 2. Auxiliary assertions.}
\vspace{0.3cm}

All the proofs of the main theorems are quite elementary. However, they rely on the following assertion that allows one to change the order of summation in multiple series.
\vspace{0.3cm}

\noindent
\textsc{Lemma 2.} \emph{Suppose that $\{\alpha_{m,n}\}$ is a sequence of complex numbers such that at least one of the following series}
\[
\sum_{m}\sum_{n}|\alpha_{m,n}|,\quad \sum_{n}\sum_{m}|\alpha_{m,n}|
\]
\emph{converges. Then the double series $\displaystyle \sum\limits_{m,n}\alpha_{m,n}$ is absolutely convergent and its sum equals}
\[
\sum_{m}\sum_{n}\alpha_{m,n} = \sum_{n}\sum_{m}\alpha_{m,n}.
\]

For the proof, see \cite[\S\S 1.62, 1.64]{Titchmarsh_1939}.

We also need several formulas for the expansions of the simplest rational functions.
\vspace{0.3cm}

\noindent
\textsc{Lemma 3.} \emph{Suppose that $k\geqslant 1$ is integer. Then, for any $s$, $0\leqslant s\leqslant 2k-1$, we have}
\begin{equation}\label{lab_02}
\frac{w^{s}}{w^{2k\mathstrut}+1} = \frac{(-1)^{s}}{2k}\sum\limits_{r=0}^{2k-1}\frac{\varepsilon_{r}^{s+1}}{w+\varepsilon_{r}}.
\end{equation}
\textsc{Proof.} First, we prove (\ref{lab_02}) for $s=0$. To do this, we define
\[
P(w) = \frac{1}{2k}\sum\limits_{r=0}^{2k-1}\frac{\varepsilon_{r}}{w-\varepsilon_{r}},\quad Q(w) = (w^{2k}+1)P(w),
\]
and fix an arbitrary $\ell$ such that $0\leqslant \ell\leqslant 2k-1$. Then
\[
Q(\varepsilon_{\ell}) = \lim_{w\to \varepsilon_{\ell}}Q(w) = \frac{\varepsilon_{\ell}}{2k}\lim_{w\to \varepsilon_{t}}\frac{w^{2k}+1}{w-\varepsilon_{\ell}}
= \frac{\varepsilon_{\ell}}{2k}\cdot 2k\varepsilon_{\ell}^{2k-1} = \varepsilon_{\ell}^{2k} = -1.
\]
Thus, the polynomial $Q(w)+1$  has at least $2k$ distinct roots, while its degree is at most $2k-1$. Hence, $Q(w)\equiv -1$ and therefore
\[
\frac{1}{w^{2k}+1} = -P(w) = -\frac{1}{2k}\sum\limits_{r=0}^{2k-1}\frac{\varepsilon_{r}}{w-\varepsilon_{r}} =
\frac{1}{2k}\sum\limits_{r=0}^{2k-1}\frac{\varepsilon_{r}}{w+\varepsilon_{r}}.
\]
Suppose now that the assertion has been proved for any $s$ such that $0\leqslant s\leqslant n$, where $n\leqslant 2k-2$. Note that for such $n$ we have
\[
\sum\limits_{r=0}^{2k-1}\varepsilon_{r}^{n+1} = \exp{\biggl(\frac{\pi i(n+1)}{2k}\biggr)}\sum\limits_{r=0}^{2k-1}\exp{\biggl(\frac{2\pi i(n+1)r}{2k}\biggr)} = 0.
\]
Hence,
\[
\frac{(-1)^{n+1}}{2k}\sum\limits_{r=0}^{2k-1}\frac{\varepsilon_{r}^{n+1}}{w} = 0.
\]
Adding this relation term by term to the identity
\[
\frac{(-1)^{n}}{2k}\sum\limits_{r=0}^{2k-1}\frac{\varepsilon_{r}^{n+1}}{w+\varepsilon_{r}} = \frac{w^{n}}{w^{2k\mathstrut}+1},
\]
we obtain
\[
\frac{w^{n}}{w^{2k}+1} = \frac{(-1)^{n+1}}{2k}\sum\limits_{r=0}^{2k-1}\varepsilon_{r}^{n+1}\biggl(\frac{1}{w}-\frac{1}{w+\varepsilon_{r}}\biggr) =
\frac{(-1)^{n+1}}{2k}\sum\limits_{r=0}^{2k-1}\frac{\varepsilon_{r}^{n+2}}{w(w+\varepsilon_{r})}.
\]
Multiplying both sides of the last equality by $w$, we arrive at the desired assertion. $\square$
\vspace{0.3cm}

\noindent
\textsc{Corollary 1.} \emph{Under the assumptions of Lemma~3, for any $m,n\geqslant 1$, we have}
\[
\frac{n^{2k-1}}{m^{2k\mathstrut}+n^{2k\mathstrut}} = \frac{1}{2k}\sum\limits_{r=0}^{k-1}\varepsilon_{r}\biggl(\frac{1}{m+n\varepsilon_{r}} - \frac{1}{m-n\varepsilon_{r}}\biggr).
\]
\textsc{Proof.} Taking $s=1$ and $w = m/n$ in (\ref{lab_02}), we obtain
\[
\frac{n^{2k-1}}{m^{2k\mathstrut}+n^{2k\mathstrut}} = \frac{1}{2k}\sum\limits_{r=0}^{2k-1}\frac{\varepsilon_{r}}{m+n\varepsilon_{r}} =
\frac{1}{2k}\sum\limits_{r=0}^{k-1}\biggl(\frac{\varepsilon_{r}}{m+n\varepsilon_{r}} + \frac{\varepsilon_{r+k}}{m+n\varepsilon_{r+k}}\biggr).
\]
It now suffices to note that $\varepsilon_{r+k} = -\varepsilon_{r}$. $\square$
\vspace{0.3cm}

\noindent
\textsc{Corollary 2.} \emph{Given the assumptions of Lemma~3, for any $m,n\geqslant 1$ and any $0\leqslant \ell\leqslant 2k-1$, we have}
\[
\frac{m^{\ell}n^{2k-1-\ell}}{m^{2k\mathstrut}+n^{2k\mathstrut}} = \frac{(-1)^{\ell+1}}{2k}\sum\limits_{r=0}^{2k-1}\varepsilon_{r}^{\ell+1}\biggl(\frac{1}{m}-\frac{1}{m+n\varepsilon_{r}}\biggr).
\]
\textsc{Proof.} Taking $s=\ell$, $w = m/n$ in (\ref{lab_02}), we obtain
\[
\frac{m^{\ell}n^{2k-1-\ell}}{m^{2k\mathstrut}+n^{2k\mathstrut}} = \frac{(-1)^{\ell}}{2k}\sum\limits_{r=0}^{2k-1}\frac{\varepsilon_{r}^{\ell+1}}{m+n\varepsilon_{r}}.
\]
Subtracting this identity from the relation
\[
\frac{(-1)^{\ell+1}}{2k}\sum\limits_{r=0}^{2k-1}\frac{\varepsilon_{r}^{\ell+1}}{m} = 0,
\]
we easily obtain the desired assertion. $\square$
\vspace{0.3cm}

\noindent
\textsc{Lemma 4.} \emph{For any $k\geqslant 1$ and $s$ such that $0\leqslant s \leqslant 2k$, one has}
\[
\frac{w^{s}}{w^{2k+1}+1} = -\frac{1}{2k+1}\sum\limits_{r=0}^{2k}\frac{\omega_{r}^{s+1}}{w-\omega_{r}}.
\]
\textsc{Corollary 1.} \emph{Under the assumptions of Lemma 4, for any $m,n\geqslant 1$ one has}
\[
\frac{n^{2k}}{m^{2k+1\mathstrut}+n^{2k+1\mathstrut}} = \frac{1}{2k+1}\sum\limits_{r=0}^{2k}\omega_{r}\biggl(\frac{1}{m} - \frac{1}{m-n\omega_{r}}\biggr).
\]
\textsc{Corollary 2.} \emph{Under the assumptions of Lemma 4, for any $m,n\geqslant 1$ and for any $k\geqslant 1$, $\ell$, $0\leqslant \ell\leqslant 2k$, one has}
\[
\frac{m^{\ell}n^{2k-\ell}}{m^{2k+1\mathstrut}+n^{2k+1\mathstrut}} = -\,\frac{1}{2k+1}\sum\limits_{r=0}^{2k}\frac{\omega_{r}^{\ell}}{m-n\omega_{r}}.
\]
The proofs of these relations follow word by word the proofs of Lemma 3 and its corollaries.

In what follows, we exploit the basic properties of the digamma function $\psi(z)$.
\vspace{0.3cm}

\noindent
\textsc{Lemma 5.} \emph{Let $z$ be a complex number such that $z\ne -m$, where $m = 0,1,2,3,\ldots.$ Then}
\[
\psi(z) = -\gamma-\frac{1}{z} + \sum\limits_{n=1}^{+\infty}\biggl(\frac{1}{n} - \frac{1}{n+z}\biggr).
\]
\textsc{Lemma 6.} \emph{For any non-integer complex $z$, one has}
\[
\psi(-z) = \psi(z) + \frac{1}{z} + \pi\ctg{(\pi z)}.
\]

\noindent
\textsc{Lemma 7.} \emph{For any complex non-integer $z$, one has}
\[
\sum\limits_{n=1}^{+\infty}\biggl(\frac{1}{n+z}-\frac{1}{n-z}\biggr) = \pi\ctg{(\pi z)} - \frac{1}{z}.
\]
These assertions are well-known; see, for example, \cite[\S 12.16]{Whittaker_Watson_1927}, \cite[\S 2.1]{Olver_1974}.
\vspace{0.3cm}

\noindent
\textsc{Lemma 8.} \emph{Suppose that $z=\varrho e^{i\varphi}$, where $\varrho>0$ and $0<|\varphi|<\pi$ are real numbers. Then the following inequality holds:}
\[
\biggl|\ctg{(\pi z)}-\frac{1}{\pi z}\biggr|< \frac{3}{2}\cdot\frac{\min{\bigl(1,\varrho\bigr)}}{|\!\sin{\varphi}|}.
\]
\textsc{Proof.} By Lemma 7, we have
\[
\biggl|\ctg{(\pi z)}-\frac{1}{\pi z}\biggr|\leqslant \frac{1}{\pi}\sum\limits_{n=1}^{+\infty}\biggl|\frac{1}{n+z}-\frac{1}{n-z}\biggr|= \frac{2|z|}{\pi}
\sum\limits_{n=1}^{+\infty}\frac{1}{|n^{2}-z^{2}|}.
\]
Furthermore,
\[
\bigl|n^{2}-z^{2}\bigr|^{2} = \bigl|n^{2}-\varrho^{2}e^{2i\varphi}\bigr|^{2} = n^{4}+\varrho^{4} - 2(n\varrho)^{2}\cos{(2\varphi)}.
\]
If $\cos{(2\varphi)}\leqslant 0$, then
\[
\bigl|n^{2}-z^{2}\bigr|^{2}\geqslant n^{4}+\varrho^{4} \geqslant \frac{1}{2}\,(n^{2}+\varrho^{2})^{2},
\]
so that
\begin{equation}\label{lab_03}
\bigl|n^{2}-z^{2}\bigr|\geqslant \frac{n^{2}+\varrho^{2}}{\sqrt{2}}.
\end{equation}
Suppose now that $\cos{(2\varphi)}>0$. Then
\[
\bigl|n^{2}-z^{2}\bigr|^{2}\geqslant n^{4}+\varrho^{4} - (n^{4}+\varrho^{4})\cos{(2\varphi)} = (n^{4}+\varrho^{4})\cdot 2\sin^{2}{\varphi}\geqslant (n^{2}+\varrho^{2})^{2}\sin^{2}{\varphi},
\]
and therefore
\begin{equation}\label{lab_04}
\bigl|n^{2}-z^{2}\bigr|\geqslant |\!\sin{\varphi}|(n^{2}+\varrho^{2}).
\end{equation}
It follows from (\ref{lab_03}) and (\ref{lab_04}) that in both cases, we have
\[
\bigl|n^{2}-z^{2}\bigr|\geqslant \frac{|\!\sin{\varphi}|}{\sqrt{2}}\,(n^{2}+\varrho^{2})
\]
Hence,
\[
\biggl|\ctg{(\pi z)}-\frac{1}{\pi z}\biggr|\leqslant \frac{2\sqrt{2}}{\pi|\!\sin{\varphi}|}\sum\limits_{n=1}^{+\infty}\frac{\varrho}{n^{2}+\varrho^{2}}.
\]
In the case $0<\varrho < 3/\pi$, the sum over $n$ is less than
\[
\varrho\sum\limits_{n=1}^{+\infty}\frac{1}{n^{2}} = \frac{\pi^{2}\varrho}{6} = \frac{\pi}{2}\cdot \frac{\pi\varrho}{3} = \frac{\pi}{2}\,\min{\biggl(1,\frac{\pi\varrho}{3}\biggr)}.
\]
If $\varrho \geqslant 3/\pi$ then this sum is at most
\[
\int_{0}^{+\infty}\frac{\varrho\,dx}{x^{2}+\varrho^{2}} = \frac{\pi}{2} = \frac{\pi}{2}\,\min{\biggl(1,\frac{\pi\varrho}{3}\biggr)}.
\]
Finally, we obtain:
\[
\biggl|\ctg{(\pi z)}-\frac{1}{\pi z}\biggr|\leqslant \frac{2\sqrt{2}}{\pi|\!\sin{\varphi}|}\cdot \frac{\pi}{2}\,\min{\biggl(1,\frac{\pi\varrho}{3}\biggr)} \leqslant \frac{\pi\sqrt{2}}{3}\cdot\frac{\min{\bigl(1,\varrho\bigr)}}{|\!\sin{\varphi}|} < \frac{3}{2}\cdot\frac{\min{\bigl(1,\varrho\bigr)}}{|\!\sin{\varphi}|}.
\]
Lemma is proved. $\square$

\textsc{Lemma 9.} \emph{For any complex $z$ such that $|\!\arg{z}|<\pi/2$ and for any integer $M\geqslant 1$, the following identity holds:}
\[
\psi(z) = \ln{z}-\frac{1}{2z}-\sum\limits_{m=1}^{M}\frac{B_{2m}}{2m}\cdot\frac{1}{z^{2m\mathstrut}} + \frac{2(-1)^{M+1}}{z^{2M\mathstrut}}\,j_{M}(z).
\]
\emph{Here}
\[
j_{M}(z) = \int_{0}^{+\infty}\frac{t^{2M+1}}{e^{2\pi t}-1}\,\frac{dt}{(t^{2}+z^{2})}
\]
\emph{and $B_{j}$ denote the Bernoulli numbers:}
\[
B_{0} = 1, \quad B_{1} = -\frac{1}{2},\quad B_{2} = \frac{1}{6},\quad B_{4} = -\frac{1}{30},\quad  B_{6} = \frac{1}{42}, \quad \ldots,
\]
\emph{with $B_{2n+1} = 0$ for all $n\geqslant 1$.}
\vspace{0.3cm}

For the proof, see \cite[\S\S 12.32, 12.33]{Whittaker_Watson_1927}
\vspace{0.3cm}

\noindent
\textsc{Lemma 10.} \emph{For any complex non-integer $z = x+iy$, one has}
\[
\ctg{(x+iy)} = \frac{\sin{(2x)}-i\sh{(2y)}}{2(\sh^{2}{(y)}+\sin^{2}{(x)})} = \frac{\sin{(2x)}-i\sh{(2y)}}{\ch{(2y)}-\cos{(2x)}}.
\]

\textbf{\S 3. Main theorems.}
\vspace{0.3cm}

Here, we prove a number of identities expressing the squares of the Riemann zeta function at integer points in terms of sums of some rapidly convergent series.
\vspace{0.3cm}

\noindent
\textsc{Theorem 1.} \emph{For any integer $k\geqslant 1$, the following relation holds:}
\[
\zeta^{2}(2k) + \zeta(4k) = \sum\limits_{n=1}^{+\infty}\frac{a_{k}(n)}{n^{4k-1\mathstrut}}.
\]
\textsc{Proof.} Fixing any $m,n\geqslant 1$ and using the Corollary 1 of Lemma 3, we obtain
\[
\frac{n^{2k-1}}{m^{2k}+n^{2k}} = \frac{1}{2k}\sum\limits_{r=0}^{k-1}\varepsilon_{r}\biggl(\frac{1}{m+n\varepsilon_{r}} - \frac{1}{m-n\varepsilon_{r}}\biggr).
\]
Hence,
\begin{equation}\label{lab_05}
\sum\limits_{m=1}^{+\infty}\frac{n^{2k-1}}{m^{2k}+n^{2k}} = \frac{1}{2k}\sum\limits_{m=1}^{+\infty}\sum\limits_{r=0}^{k-1}\varepsilon_{r}
\biggl(\frac{1}{m+n\varepsilon_{r}} - \frac{1}{m-n\varepsilon_{r}}\biggr).
\end{equation}
Further, (\ref{lab_04}) implies that
\[
\bigl|m^{2} - (n\varepsilon_{r})^{2}\bigr|\geqslant (m^{2}+n^{2})\sin{\frac{\pi}{2k}}.
\]
Hence,
\begin{equation}\label{lab_06}
\biggl|\varepsilon_{r}\biggl(\frac{1}{m+n\varepsilon_{r}} - \frac{1}{m-n\varepsilon_{r}}\biggr)\biggr| = \frac{2n}{|m^{2\mathstrut}-(n\varepsilon_{r})^{2\mathstrut}|}\leqslant \biggl(\sin{\frac{\pi}{2k\mathstrut}}\biggr)^{-1}\cdot\frac{2n}{m^{2}+n^{2\mathstrut}},.
\end{equation}
Lemma 2 together with (\ref{lab_06}) allows one to change the order of summation in the right-hand side of (\ref{lab_05}). Using the identity of Lemma 7, we find
\begin{multline}\label{lab_07}
\sum\limits_{m=1}^{+\infty}\frac{n^{2k-1}}{m^{2k\mathstrut}+n^{2k\mathstrut}} = \frac{1}{2k}\sum\limits_{r=0}^{k-1}\varepsilon_{r}\sum\limits_{m=1}^{+\infty}
\biggl(\frac{1}{m+n\varepsilon_{r}} - \frac{1}{m-n\varepsilon_{r}}\biggr) = \\ = \frac{1}{2k}\sum\limits_{r=0}^{k-1}\varepsilon_{r}\biggl(\pi\ctg{(\pi\varepsilon_{r}n)}-\frac{1}{n\varepsilon_{r}}\biggr) = \frac{1}{2}\,a_{k}(n) - \frac{1}{2n}.
\end{multline}
Further, by the definition of $a_{k}(n)$ and Lemma 8, we have
\begin{align*}
& \biggl|a_{k}(n)-\frac{1}{n}\biggr| = \biggl|\frac{\pi}{k}\sum\limits_{r=0}^{k-1}\varepsilon_{r}\biggl(\ctg{(\pi\varepsilon_{r}n)}-\frac{1}{\pi n\varepsilon_{r}}\biggr)\biggr|\leqslant\frac{\pi}{k}
\sum\limits_{r=0}^{k-1}\biggl|\ctg{(\pi\varepsilon_{r}n)}-\frac{1}{\pi n\varepsilon_{r}}\biggr|\ll_{k} 1,\\
& \bigl|a_{k}(n)\bigr|\ll_{k} 1.
\end{align*}
Multiplying both parts of (\ref{lab_07}) by $n^{-(4k-1)}$ and summing over $n\geqslant 1$, we obtain
\[
\sum\limits_{m,n=1}^{+\infty}\frac{1}{n^{2k}(m^{2k\mathstrut}+n^{2k\mathstrut})} = \frac{1}{2}\sum\limits_{n=1}^{+\infty}\frac{a_{k}(n)}{n^{4k-1}} - \frac{1}{2}\,\zeta(4k).
\]
By Lemma 1, the left-hand side equals
\[
\sum\limits_{m,n=1}^{+\infty}\frac{1}{(mn)^{2k}}\cdot\frac{m^{2k}}{m^{2k\mathstrut}+n^{2k\mathstrut}} = \frac{1}{2}\biggl(\,\sum\limits_{m=1}^{+\infty}\frac{1}{m^{2k\mathstrut}}\biggr)^{2} = \frac{1}{2}\,\zeta^{2}(2k).
\]
This completes the proof of Theorem 1. $\square$
\vspace{0.5cm}

\noindent
\textsc{Corollary.} \emph{For any $k\geqslant 1$, the following identity holds:}
\begin{equation}\label{lab_08}
\zeta^{2}(2k)+\zeta(4k) = \frac{\pi}{k}\,\frac{\zeta(4k-1)}{s(0)} + \frac{\pi}{2k}\sum\limits_{n=1}^{+\infty}\frac{\alpha_{k}(n)}{n^{4k-1\mathstrut}}.
\end{equation}
\emph{Here}
\begin{multline*}
\alpha_{k}(n) = \sum\limits_{r=0}^{k-1}\frac{c(r)\sin{(2\pi nc(r))}+s(r)\cos{(2\pi nc(r))}-s(r)e^{-2\pi ns(r)}}{\sh^{2}{(\pi ns(r))}+\sin^{2}{(\pi nc(r))}},\\
c(r) = \cos\biggl(\frac{\pi}{2k}(2r+1)\biggr),\quad s(r) = \sin\biggl(\frac{\pi}{2k}(2r+1)\biggr).
\end{multline*}
\noindent
\textsc{Remark.} Setting $k=1$ in (\ref{lab_08}) and using the relations $\zeta(2) = \pi^{2}/6$, $\zeta(4) = \pi^{4}/90$, we obtain,
after some straightforward transformations, we get the ``Cauchy-Lerch-Ramanujan formula''\footnote{This name is not generally accepted.}
\[
\zeta(3) = \frac{7\pi^{3}}{180} - 2\sum\limits_{n=1}^{+\infty}\frac{1}{n^{3\mathstrut}(e^{2\pi n}-1)}
\]
(see \cite[formula (120)]{Cauchy_1889}, \cite{Lerch_1901}, \cite[pp. 275-276, 293]{Berndt_1989} and \cite{Korolev_2026}).
\vspace{0.5cm}

\textsc{Theorem 2.} \emph{For any $k,\ell\geqslant 1$ such that $1\leqslant \ell \leqslant 2k-2$, the following relation holds:}
\[
\zeta^{2}(2k-\ell) = \sum\limits_{n=1}^{+\infty}\frac{b_{k,\ell}(n)}{n^{4k-2\ell-1}}.
\]
\textsc{Proof.} By Corollary 2 of Lemma 3, for any $m,n\geqslant 1$ we have
\[
\frac{m^{\ell}n^{2k-\ell-1}}{m^{2k\mathstrut}+n^{2k\mathstrut}} = \frac{(-1)^{\ell+1}}{2k}\sum\limits_{r=0}^{2k-1}\varepsilon_{r}^{\ell+1}\biggl(\frac{1}{m}-\frac{1}{m+n\varepsilon_{r}}\biggr).
\]
Summing both sides over $m$ and using Lemma 5, we get
\begin{multline}\label{lab_09}
\sum\limits_{m=1}^{+\infty}\frac{m^{\ell}n^{2k-\ell-1}}{m^{2k\mathstrut}+n^{2k\mathstrut}} =
\frac{(-1)^{\ell+1}}{2k}\sum\limits_{m=1}^{+\infty}\sum\limits_{r=0}^{2k-1}\varepsilon_{r}^{\ell+1}\biggl(\frac{1}{m}-\frac{1}{m+n\varepsilon_{r}}\biggr) = \\
= \frac{(-1)^{\ell+1}}{2k}\sum\limits_{r=0}^{2k-1}\varepsilon_{r}^{\ell+1}\sum\limits_{m=1}^{+\infty}\biggl(\frac{1}{m}-\frac{1}{m+n\varepsilon_{r}}\biggr) =
\frac{(-1)^{\ell+1}}{2k}\sum\limits_{r=0}^{2k-1}\varepsilon_{r}^{\ell+1}\biggl(\psi(n\varepsilon_{r})+\gamma + \frac{1}{n\varepsilon_{r}}\biggr) = \\
= \frac{(-1)^{\ell+1}}{2k}\sum\limits_{r=0}^{2k-1}\varepsilon_{r}^{\ell+1}\psi(n\varepsilon_{r}) = \frac{1}{2}\,b_{k,\ell}(n)
\end{multline}
(the change in the order of summation is justified as before).

Further, we multiply both sides of (\ref{lab_09}) by $n^{-(4k-2\ell-1)}$ and sum over $n$. Thus, by Lemma 1, we obtain
\begin{multline*}
\frac{1}{2}\sum\limits_{n=1}^{+\infty}\frac{b_{k,\ell}(n)}{n^{4k-2\ell-1}} = \sum\limits_{m,n=1}^{+\infty}\frac{1}{n^{4k-2\ell-1}}\cdot
\frac{m^{\ell}n^{2k-\ell-1}}{m^{2k\mathstrut}+n^{2k\mathstrut}} = \sum\limits_{m,n=1}^{+\infty}\frac{m^{\ell}}{n^{2k-\ell}(m^{2k\mathstrut}+n^{2k\mathstrut})} = \\
= \sum\limits_{m,n=1}^{+\infty}\frac{1}{(mn)^{2k-\ell}}\cdot\frac{m^{2k}}{(m^{2k\mathstrut}+n^{2k\mathstrut})} = \frac{1}{2}\biggl(\,\sum\limits_{m=1}^{+\infty}\frac{1}{m^{2k-\ell}}\biggr)^{2}
= \frac{1}{2}\,\zeta^{2}(2k-\ell).
\end{multline*}
Theorem is proved. $\square$
\vspace{0.3cm}

\noindent
\textsc{Corollary 1.} \emph{The following identity holds:}
\[
\zeta^{2}(3) = -\sum\limits_{n=1}^{+\infty}\frac{\beta(n)}{n^{5}},\quad \beta(n) = \IIm{\bigl(\psi(n\varepsilon_{0})+\psi(-n\varepsilon_{0})\bigr)},\quad \varepsilon_{0} = e^{\pi i/4}.
\]
\textsc{Proof.} Setting $k=2$, $\ell = 1$, we obtain
\[
\varepsilon_{0} = e^{\pi i/4},\quad \varepsilon_{1} = -\overline{\varepsilon}_{0},\quad \varepsilon_{2} = -\varepsilon_{0},\quad
\varepsilon_{3} = \overline{\varepsilon}_{0}.
\]
Further, we have $\varepsilon_{r}^{2} = (-1)^{r}\cdot i$ for $r=0,1,2$. Hence,
\begin{multline*}
b_{2}(n) = \frac{1}{2}\sum\limits_{r=0}^{3}\varepsilon_{r}^{2}\psi(n\varepsilon_{r}) = \frac{i}{2}\sum\limits_{r=0}^{3}(-1)^{r}\psi(n\varepsilon_{r}) = \\
= \frac{i}{2}\bigl(\psi(n\varepsilon_{0})-\psi(n\varepsilon_{1})+\psi(n\varepsilon_{2})-\psi(n\varepsilon_{3})\bigr) =
\frac{i}{2}\bigl(\psi(n\varepsilon_{0})-\psi(-n\overline{\varepsilon}_{0})+\psi(-n\varepsilon_{0})-\psi(n\overline{\varepsilon}_{0})\bigr) = \\
= -\IIm{\bigl(\psi(n\varepsilon_{0})+\psi(-n\varepsilon_{0})\bigr)}.\quad \square
\end{multline*}

\noindent
\textsc{Corollary 2.} \emph{For any fixed $m\geqslant 0$, the following identity holds:}
\[
\zeta^{2}(3) = \frac{\pi}{2}\,\zeta(5) - \sum\limits_{r=0}^{m}\frac{(-1)^{r}B_{4r+2}}{2r+1}\,\zeta(4r+7) + \mathfrak{S}_{0} + (-1)^{m}\int_{0}^{+\infty}\frac{G_{m}(t)dt}{e^{2\pi t}-1}.
\]
\emph{Here}
\[
\mathfrak{S}_{0} = \pi\sum\limits_{n=1}^{+\infty}\frac{1}{n^{5}}\cdot \frac{\cos{(\pi n\sqrt{2})}-e^{-\pi n\sqrt{2}}}{\ch{(\pi n\sqrt{2})}-\cos{(\pi n\sqrt{2})}} = -0.0204388172\ldots,
\quad G_{m}(t) = 4\sum\limits_{n=1}^{\infty}\frac{(t/n)^{4m+5}}{n^{2\mathstrut}(n^{4\mathstrut}+t^{4\mathstrut})}.
\]
\textsc{Proof.} By Lemma 6,
\[
\psi(n\varepsilon_{0})+\psi(-n\varepsilon_{0}) = 2\psi(n\varepsilon_{0})+\frac{1}{n\varepsilon_{0}} + \pi\ctg{(\pi\varepsilon_{0}n)}.
\]
Further, Lemma 10 implies that
\[
\ctg{(\pi\varepsilon_{0}n)} = \frac{\sin{(\pi n\sqrt{2})}-i\sh{(\pi n\sqrt{2})}}{\ch{(\pi n\sqrt{2})}-\cos{(\pi n\sqrt{2})}}.
\]
Hence,
\begin{multline*}
\beta(n) = 2\IIm{\psi(n\varepsilon_{0})} - \frac{1}{n\sqrt{2}} - \frac{\pi\sh{(\pi n\sqrt{2})}}{\ch{(\pi n\sqrt{2})}-\cos{(\pi n\sqrt{2})}} = \\
= 2\IIm{\psi(n\varepsilon_{0})} - \frac{1}{n\sqrt{2}} - \pi\biggl(1 + \frac{\cos{(\pi n\sqrt{2})}-e^{-\pi n\sqrt{2}}}{\ch{(\pi n\sqrt{2})}-\cos{(\pi n\sqrt{2})}}\biggr).
\end{multline*}
Taking $z = n\varepsilon_{0}$, $M = 2m+1$ in Lemma 9, we find
\[
\psi(n\varepsilon_{0}) = \ln{n} + \frac{\pi i}{4} + \frac{i-1}{2n\sqrt{2}} - \sum\limits_{s=1}^{2m+1}\frac{B_{2s}}{2s}\cdot \frac{i^{-s}}{n^{2s}} + (-1)^{m+1}\cdot\frac{2i}{n^{4m+2\mathstrut}}
\int_{0}^{+\infty}\frac{t^{4m+3}}{e^{2\pi t}-1}\cdot\frac{dt}{t^{2\mathstrut}+in^{2\mathstrut}}
\]
and therefore
\[
\IIm \psi(n\varepsilon_{0}) = \frac{\pi}{4} + \frac{1}{2n\sqrt{2}} -  \sum\limits_{s=1}^{2m+1}\frac{B_{2s}}{2s}\cdot \frac{\IIm{(i^{-s})}}{n^{2s}} +
(-1)^{m+1}\cdot\frac{2}{n^{4m+2\mathstrut}}
\int_{0}^{+\infty}\frac{t^{4m+5}}{e^{2\pi t}-1}\cdot\frac{dt}{t^{4\mathstrut}+n^{4\mathstrut}}.
\]
The non-zero contribution to the sum over $s$ is given by odd $s = 2r+1$ only. For such $s$, we have $\IIm{(i^{-s})} = (-1)^{r+1}$. Thus we obtain:
\begin{multline*}
\beta(n) = -\frac{\pi}{2} + \sum\limits_{r=0}^{m}\frac{(-1)^{r}B_{4r+2}}{2r+1}\cdot \frac{1}{n^{4r+2}} - \pi\,\frac{\cos{(\pi n\sqrt{2})}-e^{-\pi n\sqrt{2}}}{\ch{(\pi n\sqrt{2})}-\cos{(\pi n\sqrt{2})}} + \\
+ \frac{4(-1)^{m+1}}{n^{4m+2\mathstrut}}\int_{0}^{+\infty}\frac{t^{4m+5}}{e^{2\pi t}-1}\cdot\frac{dt}{t^{4\mathstrut}+n^{4\mathstrut}}.
\end{multline*}
Now the assertion follows. $\square$
\vspace{0.5cm}

\noindent
\textsc{Theorem 3.} \emph{For any $k\geqslant 1$, the following identity holds:}
\[
\frac{1}{2}\,\zeta^{2}(2k+1) + \zeta(4k+2) = \sum\limits_{n=1}^{+\infty}\frac{c_{k}(n)}{n^{4k+1}}.
\]
\textsc{Proof.} By Corollary 1 of Lemma 4, for any $m,n\geqslant 1$, one has:
\begin{equation}\label{lab_10}
\frac{n^{2k}}{n^{2k+1\mathstrut}+m^{2k+1\mathstrut}} = \frac{1}{2k+1}\sum\limits_{r=0}^{2k}\omega_{r}\biggl(\frac{1}{m} - \frac{1}{m-n\omega_{r}}\biggr).
\end{equation}
Summing both sides of (\ref{lab_10}) over $m\geqslant 1$, and using Lemma 5, we find
\begin{multline}\label{lab_11}
\sum\limits_{m=1}^{+\infty}\frac{n^{2k}}{n^{2k+1\mathstrut}+m^{2k+1\mathstrut}} = \frac{1}{2k+1}\sum\limits_{r=0}^{2k}\omega_{r}\biggl(\psi(-n\omega_{r})+\gamma - \frac{1}{n\omega_{r}}\biggr) = \\ = \frac{1}{2k+1}\sum\limits_{r=0}^{2k}\omega_{r}\psi(-n\omega_{r}) - \frac{1}{n} = c_{k}(n) - \frac{1}{n}.
\end{multline}
The change in the order of summation is justified as before, but here we use the inequality
\[
|m-n\omega_{r}|\geqslant (m+n)\sin{\frac{\pi}{2(2k+1)}}.
\]
Further,
\begin{multline*}
\biggl|c_{k}(n) - \frac{1}{n}\biggr| \leqslant \frac{1}{2k+1}\sum\limits_{r=0}^{2k}\biggl|\psi(-n\omega_{r})+\gamma - \frac{1}{n\omega_{r}}\biggr| =
\frac{1}{2k+1}\sum\limits_{r=0}^{2k}\biggl|\sum\limits_{m=1}^{+\infty}\biggl(\frac{1}{m} - \frac{1}{m-n\varepsilon_{r}}\biggr)\biggr| \leqslant \\
\leqslant \frac{1}{2k+1}\sum\limits_{r=0}^{2k}\sum\limits_{m=1}^{+\infty}\frac{n}{m|m-n\omega_{r}|}\leqslant \biggl((2k+1)\sin{\frac{\pi}{2(2k+1)}}\biggr)^{-1}
\sum\limits_{r=0}^{2k}\sum\limits_{m=1}^{+\infty}\frac{n}{m(m+n)} = \\
= \biggl(\sin{\frac{\pi}{2(2k+1)}}\biggr)^{-1}\sum\limits_{m=1}^{+\infty}\biggl(\frac{1}{m}-\frac{1}{m+n}\biggr) =
\biggl(\sin{\frac{\pi}{2(2k+1)}}\biggr)^{-1}H_{n},
\end{multline*}
where
\[
H_{n} = 1+\frac{1}{2} + \frac{1}{3} + \ldots + \frac{1}{n}\leqslant \ln{n}+1.
\]
Hence,
\begin{equation}\label{lab_12}
|c_{k}(n)|\leqslant \ln{n}+1+\frac{1}{n}\leqslant \ln{n}+2.
\end{equation}
Multiplying both sides of (\ref{lab_11}) by $n^{-(4k+1)}$ and summing over $n$, we obtain
\begin{equation}\label{lab_13}
\sum\limits_{n=1}^{+\infty}\frac{1}{n^{4k+1}}\biggl(c_{k}(n)-\frac{1}{n}\biggr) = \sum\limits_{m,n=1}^{+\infty}\frac{1}{n^{4k+1}}\cdot \frac{n^{2k}}{m^{2k+1\mathstrut}+n^{2k+1\mathstrut}}.
\end{equation}
In view of the estimate (\ref{lab_12}), the left-hand side of (\ref{lab_13}) is equal to
\[
\sum\limits_{n=1}^{+\infty}\frac{c_{k}(n)}{n^{4k+1\mathstrut}} - \zeta(4k+2).
\]
By Lemma 1, the right-hand side of (\ref{lab_13}) coincides with
\[
\sum\limits_{m,n=1}^{+\infty}\frac{m^{2k+1}}{(mn)^{2k+1\mathstrut}(m^{2k+1\mathstrut}+n^{2k+1\mathstrut})} = \frac{1}{2}\biggl(\,\sum\limits_{m=1}^{+\infty}\frac{1}{m^{2k+1\mathstrut}}\biggr)^{2}
=\frac{1}{2}\,\zeta^{2}(2k+1).
\]
This establishes the theorem. $\square$
\vspace{0.5cm}

\noindent
\textsc{Corollary 1.} \emph{The following identity holds:}
\[
\zeta^{2}(3) = \frac{2\pi}{\sqrt{3}}\,\zeta(5) - \frac{2}{3}\,\zeta(6) + \mathfrak{S} + 2\sum\limits_{n=1}^{+\infty}\frac{\alpha(n)}{n^{5}}.
\]
\emph{Here}
\begin{equation*}
\mathfrak{S} = \frac{2\pi}{\sqrt{3}}\sum\limits_{n=1}^{+\infty}\frac{(-1)^{n}}{n^{5}}\cdot \frac{e^{-\pi n\sqrt{3}/2}}{\varphi_{n}(\pi n\sqrt{3}/2)} = -0.0312999121\ldots,
\quad
\varphi_{n}(x) =
\begin{cases}
\sh{x} & \text{\emph{for even}} \;\;n, \\
\ch{x} & \text{\emph{for odd}} \;\;n, \\
\end{cases}
\end{equation*}
\[
\alpha(n) = \frac{2}{3}\,\RRe{\bigl(\omega_{0}\psi(n\omega_{0})\bigr)} - \frac{1}{3}\,\psi(n),\quad \omega_{0} = e^{\pi i/3}.
\]
\textsc{Proof.} Taking $k=1$, we obtain
\[
\omega_{r} = e^{\pi i(2r+1)/3}\quad\text{and therefore}\quad \omega_{0} = e^{\pi i/3},\quad \omega_{1} = -1,\quad \omega_{2} = \overline{\omega}_{0} = e^{-\pi i/3}.
\]
Hence,
\[
c_{1}(n) = \frac{1}{3}\sum\limits_{r=0}^{2}\omega_{r}\psi(-n\omega_{r}) = \frac{2}{3}\,\RRe{\bigl(\omega_{0}\psi(-n\omega_{0})\bigr)} - \frac{1}{3}\,\psi(n).
\]
By Lemma 6, we have
\[
\omega_{0}\psi(-n\omega_{0}) = \omega_{0}\psi(n\omega_{0}) + \frac{1}{n} + \pi\omega_{0}\ctg{(\pi\omega_{0}n)},
\]
so that
\[
c_{1}(n) = \alpha(n) + \frac{2}{3n} + \frac{2}{3}\,\RRe{\bigl(\pi\omega_{0}\ctg{(\pi\omega_{0}n)}\bigr)}.
\]
Next, Lemma 10 implies that
\[
\ctg{(\pi\omega_{0}n)} = \ctg{\biggl(\frac{\pi n}{2} + \frac{\pi in}{2}\sqrt{3}\biggr)} = -\,\frac{i\sh{(\pi n\sqrt{3})}}{\ch(\pi n\sqrt{3})-\cos(\pi n)}
\]
and therefore,
\begin{multline*}
\RRe{\bigl(\pi\omega_{0}\ctg{(\pi\omega_{0}n)}\bigr)} = \frac{\pi\sqrt{3}}{2}\cdot \frac{\sh{(\pi n\sqrt{3})}}{\ch(\pi n\sqrt{3})-\cos(\pi n)} = \\
= \frac{\pi\sqrt{3}}{2}\cdot\biggl(1 + \frac{\cos{(\pi n)}-e^{-\pi n\sqrt{3}}}{\ch(\pi n\sqrt{3})-\cos(\pi n)}\biggr).
\end{multline*}
Now it is not difficult to verify that the second term in brackets coincides with
\[
\frac{(-1)^{n}e^{-\pi n\sqrt{3}/2}}{\varphi_{n}\bigl(\pi n\sqrt{3}/2\bigr)}.
\]
Thus,
\[
c_{1}(n) = \alpha(n) + \frac{2}{3n} + \frac{\pi}{\sqrt{3}}\biggl(1 + \frac{(-1)^{n}e^{-\pi n\sqrt{3}/2}}{\varphi_{n}\bigl(\pi n\sqrt{3}/2\bigr)}\biggr).
\]
The desired assertion now follows. $\square$
\vspace{0.5cm}

\noindent
\textsc{Corollary 2.} \emph{For any fixed $m\geqslant 0$, the following identity holds:}
\[
\zeta^{2}(3) + \zeta(6) = \frac{4\pi}{3\sqrt{3}}\,\zeta(5) + \sum\limits_{r=0}^{m}\frac{B_{6r+4}}{3r+2}\,\zeta(6r+9) + \mathfrak{S} + (-1)^{m}\int_{0}^{+\infty}\frac{F_{m}(t)dt}{e^{2\pi t}-1},
\]
\emph{where $\mathfrak{S}$ is defined in Corollary 1 and}
\[
F_{m}(t) = 4\sum\limits_{n=1}^{\infty}\frac{(t/n)^{6m+9}}{n^{6\mathstrut}+t^{6\mathstrut}}.
\]
\textsc{Proof.} Setting $M = 3m+2$ in Lemma 9, after straightforward calculations we obtain
\begin{multline*}
\RRe{\bigl(\omega_{0}\psi(n\omega_{0})\bigr)} = \frac{1}{2}\,\ln{n} - \frac{\pi}{2\sqrt{3}} - \frac{1}{2n} - \sum\limits_{s=1}^{3m+2}\frac{B_{2s}}{2s}\cdot \frac{\RRe{\bigl(\omega_{0}^{1-2s}\bigr)}}{n^{2s\mathstrut}} + \\ + \frac{2(-1)^{m}}{n^{6m+4}}\int_{0}^{+\infty}\RRe{\biggl(\frac{1}{t^{2}+n^{2}\omega_{0}^{2}}\biggr)}\frac{t^{6m+5}\,dt}{e^{2\pi t\mathstrut}-1},\\
\psi(n) = \ln{n} - \frac{1}{2n} - \sum\limits_{s=1}^{3m+2}\frac{B_{2s}}{2s}\cdot\frac{1}{n^{2s}} +
\frac{2(-1)^{m}}{n^{6m+4}}\int_{0}^{+\infty}\frac{t^{6m+5}}{e^{2\pi t\mathstrut}-1}\,\frac{dt}{t^{2}+n^{2}}.
\end{multline*}
Thus we find that
\[
\alpha(n) = -\,\frac{\pi}{3\sqrt{3}} - \frac{1}{6n} + \frac{1}{6}\sum\limits_{s=1}^{3m+2}\frac{B_{2s}}{s n^{2s}}\bigl(1-2\RRe{(\omega_{0}^{1-2s})}\bigr) + \frac{2}{3}\frac{(-1)^{m}}{n^{6m+4}}
\int_{0}^{+\infty}\frac{t^{6m+5}}{e^{2\pi t\mathstrut}-1}\,g_{n}(t)dt,
\]
where
\[
g_{n}(t) = \frac{1}{t^{2}+n^{2}} + 2\RRe{\biggl(\frac{1}{t^{2}+n^{2}\omega_{0}^{2}}\biggr)}
\]
Note that
\begin{equation*}
1-2\RRe{\bigl(\omega_{0}^{1-2s}\bigr)} = 1-2\cos{\frac{\pi}{3}(2s-1)} =
\begin{cases}
0, & \text{ if } s\not\equiv 2\pmod{3},\\
3, & \text{ if } s\equiv 2\pmod{3}.
\end{cases}
\end{equation*}
Hence, the non-zero contribution to the sum over $s$ is given by the terms corresponding to $s = 3r+2$ for $0\leqslant r\leqslant m$. Further,
\[
g_{n}(t) = \frac{1}{t^{2\mathstrut}+n^{2\mathstrut}} + \frac{2(t^{2}-n^{2}/2)}{t^{4\mathstrut}-(tn)^{2\mathstrut}+n^{4\mathstrut}} = \frac{3t^{4}}{t^{6\mathstrut}+n^{6\mathstrut}}.
\]
Thus we find
\[
\alpha(n) = -\frac{\pi}{3\sqrt{3}} - \frac{1}{6n} + \frac{1}{2}\sum\limits_{r=0}^{m}\frac{B_{6r+4}}{3r+2}\,\frac{1}{n^{6r+4}} + \frac{2(-1)^{m}}{n^{6m+4}}
\int_{0}^{+\infty}\frac{t^{6m+9}}{e^{2\pi t\mathstrut}-1}\,\frac{dt}{t^{6\mathstrut}+n^{6\mathstrut}}.
\]
Plugging this into the formula for $\zeta^{2}(3)$ from Corollary 1, we obtain:
\begin{multline*}
\zeta^{2}(3) = \frac{2\pi}{\sqrt{3}}\,\zeta(5) - \frac{2}{3}\,\zeta(6) + \mathfrak{S} + \\
+ \sum\limits_{n=1}^{+\infty}\frac{1}{n^{5}}\biggl(-\frac{2\pi}{3\sqrt{3}} - \frac{1}{3n} + \sum\limits_{r=0}^{m}\frac{B_{6r+4}}{3r+2}\,\frac{1}{n^{6r+4}} +
\frac{4(-1)^{m}}{n^{6m+4}}
\int_{0}^{+\infty}\frac{t^{6m+9}}{e^{2\pi t\mathstrut}-1}\,\frac{dt}{t^{6\mathstrut}+n^{6\mathstrut}}\biggr) = \\
= \frac{2\pi}{\sqrt{3}}\,\zeta(5) - \frac{2}{3}\,\zeta(6) + \mathfrak{S} - \frac{2\pi}{3\sqrt{3}}\,\zeta(5) - \frac{1}{3}\,\zeta(6) + \sum\limits_{r=0}^{m}\frac{B_{6r+4}}{3r+2}\,\zeta(6r+9) + \\
+ (-1)^{m}\int_{0}^{+\infty}\frac{4}{e^{2\pi t\mathstrut}-1}\sum\limits_{n=1}^{+\infty}\frac{(t/n)^{6}}{t^{6\mathstrut}+n^{6\mathstrut}}\,dt = \\
= \frac{4\pi}{3\sqrt{3}}\,\zeta(5) - \zeta(6) + \sum\limits_{r=0}^{m}\frac{B_{6r+4}}{3r+2}\,\zeta(6r+9) +
(-1)^{m}\int_{0}^{+\infty}\frac{F_{m}(t)\,dt}{e^{2\pi t\mathstrut}-1}.
\end{multline*}
Corollary is proved. $\square$
\vspace{0.5cm}

Setting $m = 0,1,2$, we obtain the following identities:
\begin{align*}
\zeta^{2}(3) + \zeta(6) = & \frac{4\pi}{3\sqrt{3}}\,\zeta(5) - \frac{1}{60}\,\zeta(9) + \mathfrak{S} + \int_{0}^{+\infty}\frac{F_{0}(t)\,dt}{e^{2\pi t\mathstrut}-1}, \\
& F_{0}(t) = 4\sum\limits_{n=1}^{+\infty}\biggl(\frac{t}{n}\biggr)^{9}\,\frac{1}{t^{6\mathstrut}+n^{6\mathstrut}},
\end{align*}
\begin{align*}
\zeta^{2}(3) + \zeta(6) = & \frac{4\pi}{3\sqrt{3}}\,\zeta(5) - \frac{1}{60}\,\zeta(9) + \frac{1}{66}\,\zeta(15) + \mathfrak{S} - \int_{0}^{+\infty}\frac{F_{1}(t)\,dt}{e^{2\pi t\mathstrut}-1}, \\
& F_{1}(t) = 4\sum\limits_{n=1}^{+\infty}\biggl(\frac{t}{n}\biggr)^{15}\,\frac{1}{t^{6\mathstrut}+n^{6\mathstrut}},\\
\zeta^{2}(3) + \zeta(6) = & \frac{4\pi}{3\sqrt{3}}\,\zeta(5) - \frac{1}{60}\,\zeta(9) + \frac{1}{66}\,\zeta(15) - \frac{3617}{4080}\,\zeta(21) + \mathfrak{S} - \int_{0}^{+\infty}\frac{F_{2}(t)\,dt}{e^{2\pi t\mathstrut}-1}, \\
& F_{2}(t) = 4\sum\limits_{n=1}^{+\infty}\biggl(\frac{t}{n}\biggr)^{21}\,\frac{1}{t^{6\mathstrut}+n^{6\mathstrut}}.
\end{align*}
\vspace{0.5cm}

\textbf{\S 4. Identities involving double series.}
\vspace{0.3cm}

We first establish a general identity.

\vspace{0.3cm}
\textsc{Theorem 4.} \emph{Let $g$ be an arithmetic function and let $f$ be its Dirichlet convolution, that is,}
\[
f(n) = \sum\limits_{d|n}g(d),\quad \textit{for any } n\geqslant 1.
\]
\emph{Suppose also that $f(n), g(n)\ll _{\alpha} n^{\alpha}$ for some fixed $\alpha$ such that $0<\alpha<1$. Finally, for $\RRe{z} > 1+\alpha$, let}
\[
L(z;f) = \sum\limits_{n=1}^{+\infty}\frac{f(n)}{n^{z}}.
\]
\emph{Then, for any $k\geqslant 1$, the following identity holds:}
\[
L^{2}(2k;f) = \sum\limits_{m=1}^{+\infty}\frac{g(m)}{m}\biggl(\,\sum\limits_{n=1}^{+\infty}\frac{f(n)}{n^{4k-1\mathstrut}}\,a_{k}\Bigl(\frac{n}{m}\Bigr)-mL(4k;f)\biggr).
\]
\textsc{Proof.} By the Corollary 1 of Lemma 3, for any $n,q\geqslant 1$, one has
\begin{multline*}
\frac{n^{2k-1}}{q^{2k\mathstrut}+n^{2k\mathstrut}} = \frac{1}{2k}\sum\limits_{r=0}^{k-1}\varepsilon_{r}\biggl(\frac{1}{q+n\varepsilon_{r}} - \frac{1}{q-n\varepsilon_{r}}\biggr) =
\frac{1}{2k}\sum\limits_{r=0}^{k-1}\biggl(\frac{1}{n+q\overline{\varepsilon}_{r}} + \frac{1}{n-q\overline{\varepsilon}_{r}}\biggr) = \\
= \frac{1}{2k}\sum\limits_{r=0}^{k-1}\biggl(\frac{1}{n+q\varepsilon_{r}} + \frac{1}{n-q\varepsilon_{r}}\biggr).
\end{multline*}
Multiplying both sides of this identity by $f(q)$ and summing over $q\geqslant 1$, we obtain
\begin{equation}\label{lab_14}
\sum\limits_{q=1}^{+\infty}\frac{n^{2k-1}f(q)}{q^{2k\mathstrut}+n^{2k\mathstrut}} = \frac{1}{2k}\sum\limits_{q=1}^{+\infty}f(q)\sum\limits_{r=0}^{k-1}\biggl(\frac{1}{n+q\varepsilon_{r}} + \frac{1}{n-q\varepsilon_{r}}\biggr).
\end{equation}
Since
\begin{equation}\label{lab_15}
\biggl|\frac{1}{n+q\varepsilon_{r}} + \frac{1}{n-q\varepsilon_{r}}\biggr| = \frac{2n}{|n^{2\mathstrut}-(q\varepsilon_{r})^{2\mathstrut}|}\leqslant \biggl(\sin{\frac{\pi}{2k}}\biggr)^{-1}\frac{2n}{n^{2}+q^{2}},
\end{equation}
the series
\[
\sum\limits_{q=1}^{+\infty}f(q)\biggl(\frac{1}{n+q\varepsilon_{r}} + \frac{1}{n-q\varepsilon_{r}}\biggr)
\]
converges absolutely for any $n\geqslant 1$, and for any $r$ such that $0\leqslant r\leqslant k-1$. Hence, Lemma 2 implies that
\begin{equation}\label{lab_16}
\sum\limits_{q=1}^{+\infty}\frac{n^{2k-1}f(q)}{q^{2k\mathstrut}+n^{2k\mathstrut}} = \frac{1}{2k}\sum\limits_{r=0}^{k-1}\sum\limits_{q=1}^{+\infty}f(q)\biggl(\frac{1}{n+q\varepsilon_{r}} + \frac{1}{n-q\varepsilon_{r}}\biggr).
\end{equation}
Now let us choose a sufficiently large $M$ and split each sum over $q$ into parts corresponding to the ranges $1\leqslant q\leqslant M$ and $q>M$.
The contribution given by $q>M$ is estimated by
\[
\sum\limits_{q>M}\frac{2n|f(q)|}{n^{2}+q^{2}}\,\biggl(\sin{\frac{\pi}{2k}}\biggr)^{-1}\ll_{k} n\sum\limits_{q>M}\frac{|f(q)|}{q^{2}}\ll_{k,\alpha}\frac{n}{M^{1-\alpha}}.
\]
Next, the sum over $1\leqslant q\leqslant M$ is transformed as follows:
\begin{multline}\label{lab_17}
\sum\limits_{1\leqslant q\leqslant M}\biggl(\,\sum\limits_{m|q}g(m)\biggr)\biggl(\frac{1}{n+q\varepsilon_{r}} + \frac{1}{n-q\varepsilon_{r}}\biggr) = \\
 = \sum\limits_{1\leqslant \ell m\leqslant M}g(m)\biggl(\frac{1}{n+\ell m\varepsilon_{r}} + \frac{1}{n-\ell m\varepsilon_{r}}\biggr) = \\
 =\sum\limits_{1\leqslant m\leqslant M}g(m)\sum\limits_{1\leqslant \ell \leqslant M/m}\biggl(\frac{1}{n+\ell m\varepsilon_{r}} + \frac{1}{n-\ell m\varepsilon_{r}}\biggr).
\end{multline}
It follows from (\ref{lab_15}) that
\[
\sum\limits_{\ell > M/m}\biggl(\frac{1}{n+\ell m\varepsilon_{r}} + \frac{1}{n-\ell m\varepsilon_{r}}\biggr)\ll_{k}
\sum\limits_{\ell > M/m}\frac{n}{n^{2\mathstrut}+(\ell m)^{2\mathstrut}}\ll_{k} \frac{n}{m^{2\mathstrut}}\sum\limits_{\ell > M/m}\frac{1}{\ell^{2}}\ll_{k}\frac{n}{Mm}.
\]
Therefore, the error introduced by replacing the sum over $\ell$ in (\ref{lab_17}) with infinite is at most
\[
\frac{n}{M}\sum\limits_{m\leqslant M}\frac{|g(m)|}{m}\ll_{k,\alpha} \frac{n}{M}\sum\limits_{m\leqslant M}\frac{1}{m^{1-\alpha}}\ll_{k,\alpha} \frac{n}{M^{1-\alpha}}
\]
Thus, the sum (\ref{lab_17}) coincides with
\begin{equation}\label{lab_18}
\sum\limits_{m\leqslant M}g(m)\sum\limits_{\ell =1}^{+\infty}\biggl(\frac{1}{n+\ell m\varepsilon_{r}} + \frac{1}{n-\ell m\varepsilon_{r}}\biggr)+ O\biggl(\frac{n}{M^{1-\alpha}}\biggr).
\end{equation}
In view of (\ref{lab_15}), for any $m$, the sum over $\ell$ is bounded by
\[
\biggl(\sin{\frac{\pi}{2k}}\biggr)^{-1}\sum\limits_{\ell =1}^{+\infty}\frac{2n}{n^{2}+(\ell m)^{2}}\ll_{k} \sum\limits_{\ell = 1}^{+\infty}\frac{n}{n^{2\mathstrut}+(\ell m)^{2\mathstrut}}
\ll_{k}\frac{n}{m^{2\mathstrut}}.
\]
Hence, the error that occurs when replacing the sum over $m$ in (\ref{lab_18}) with an infinite one is estimated as
\[
n\sum\limits_{m>M}\frac{|g(m)|}{m^{2}}\ll_{k,\alpha} \frac{n}{M^{1-\alpha}}.
\]
Therefore, the right-hand side of (\ref{lab_17}) coincides with
\[
\sum\limits_{m=1}^{+\infty}g(m)\sum\limits_{\ell = 1}^{+\infty}\biggl(\frac{1}{n+\ell m\varepsilon_{r}} + \frac{1}{n-\ell m\varepsilon_{r}}\biggr)+ O\biggl(\frac{n}{M^{1-\alpha}}\biggr).
\]
As $M\to +\infty$, we obtain
\begin{equation}\label{lab_19}
\sum\limits_{q=1}^{+\infty}f(q)\biggl(\frac{1}{n+q\varepsilon_{r}} + \frac{1}{n-q\varepsilon_{r}}\biggr) = \sum\limits_{m=1}^{+\infty}g(m)
\sum\limits_{\ell = 1}^{+\infty}\biggl(\frac{1}{n+\ell m\varepsilon_{r}} + \frac{1}{n-\ell m\varepsilon_{r}}\biggr).
\end{equation}
By Lemma 7, for any $m,n,r$, the sum over $\ell$ in (\ref{lab_19}) reduces to
\begin{multline*}
\frac{1}{m\varepsilon_{r}}\sum\limits_{\ell = 1}^{+\infty}\biggl(\frac{1}{\ell + n\overline{\varepsilon}_{r}/m} - \frac{1}{\ell - n\overline{\varepsilon}_{r}/m}\biggr) =
\frac{\overline{\varepsilon}_{r}}{m}\biggl(\pi\ctg{\frac{\pi \overline{\varepsilon}_{r}n}{m}} - \frac{m}{n\overline{\varepsilon}_{r}}\biggr) = \\
= \frac{1}{m}\biggl(\pi\overline{\varepsilon}_{r}\ctg{\frac{\pi \overline{\varepsilon}_{r}n}{m}} - \frac{m}{n}\biggr).
\end{multline*}
In view of Lemma 8, for any $m,n,r$, this quantity is estimated as
\[
\ll \biggl(\sin{\frac{\pi}{2k+1}}\biggr)^{-1}\frac{1}{m}\,\min{\biggl(1,\frac{n}{m}\biggr)}\ll_{k} \frac{\min{(m,n)}}{m^{2}}.
\]
Thus, the series
\begin{equation}\label{lab_20}
\sum\limits_{m=1}^{+\infty}g(m)\sum\limits_{\ell = 1}^{+\infty}\biggl(\frac{1}{n+\ell m\varepsilon_{r}} + \frac{1}{n-\ell m\varepsilon_{r}}\biggr) =
\sum\limits_{m=1}^{+\infty}\frac{g(m)}{m}\biggl(\pi\overline{\varepsilon}_{r}\ctg{\frac{\pi \overline{\varepsilon}_{r}n}{m}} - \frac{m}{n}\biggr)
\end{equation}
converges absolutely. By Lemma 2, after changing the order of summation in (\ref{lab_16}), (\ref{lab_19}), and (\ref{lab_20}), we obtain
\begin{equation}\label{lab_21}
\sum\limits_{q=1}^{+\infty}\frac{n^{2k-1}f(q)}{q^{2k\mathstrut}+n^{2k\mathstrut}} = \frac{1}{2}\sum\limits_{m=1}^{+\infty}\frac{g(m)}{m}\,\frac{1}{k}\sum\limits_{r=0}^{k-1}
\biggl(\pi\overline{\varepsilon}_{r}\ctg{\frac{\pi \overline{\varepsilon}_{r}n}{m}} - \frac{m}{n}\biggr).
\end{equation}
Since
\[
\biggl|\pi\overline{\varepsilon}_{r}\ctg{\frac{\pi \overline{\varepsilon}_{r}n}{m}} - \frac{m}{n}\biggr|\ll_{k} \min{\biggl(1,\frac{n}{m}\biggr)}\ll_{k}  \frac{\min{(m,n)}}{m},
\]
it follows that
\begin{equation}\label{lab_22}
\frac{1}{k}\sum\limits_{r=0}^{k-1}\biggl(\pi\overline{\varepsilon}_{r}\ctg{\frac{\pi \overline{\varepsilon}_{r}n}{m}} - \frac{m}{n}\biggr) = \overline{a}_{k}\biggl(\frac{n}{m}\biggr) - \frac{m}{n} = a_{k}\biggl(\frac{n}{m}\biggr) - \frac{m}{n} \ll \frac{\min{(m,n)}}{m}
\end{equation}
and hence
\begin{equation}\label{lab_23}
\biggl|a_{k}\biggl(\frac{n}{m}\biggr)\biggr|\ll \frac{m}{n} +  \frac{\min{(m,n)}}{m}.
\end{equation}
Multiplying both parts of (\ref{lab_21}) by $f(n)\cdot n^{-(4k-1)}$ and summing over $n\geqslant 1$, we find
\begin{equation}\label{lab_24}
\sum\limits_{n=1}^{+\infty}\frac{f(n)}{n^{2k}}\sum\limits_{q=1}^{+\infty}\frac{f(q)}{q^{2k\mathstrut}+n^{2k\mathstrut}}  =
\frac{1}{2}\sum\limits_{n=1}^{+\infty}\frac{f(n)}{n^{4k-1}}\sum\limits_{m=1}^{+\infty}\frac{g(m)}{m}\biggl(a_{k}\biggl(\frac{n}{m}\biggr) - \frac{m}{n}\biggr).
\end{equation}
By Lemma 1, the left-hand side of (\ref{lab_24}) is equal to
\[
\sum\limits_{n,q=1}^{+\infty}\frac{f(n)f(q)}{(nq)^{2k\mathstrut}}\cdot\frac{q^{2k}}{q^{2k\mathstrut}+n^{2k\mathstrut}} =
\frac{1}{2}\biggl(\,\sum\limits_{q=1}^{+\infty}\frac{f(q)}{q^{2k\mathstrut}}\biggr)^{2} = \frac{1}{2}\,L^{2}(2k;f).
\]
In view of (\ref{lab_22}), the double series in the right-hand side of (\ref{lab_24}) converges absolutely. By Lemma 2, we obtain
\[
\frac{1}{2}\sum\limits_{n=1}^{+\infty}\frac{f(n)}{n^{4k-1}}\sum\limits_{m=1}^{+\infty}\frac{g(m)}{m}\biggl(a_{k}\biggl(\frac{n}{m}\biggr) - \frac{m}{n}\biggr) = \frac{1}{2}\sum\limits_{m=1}^{+\infty}\frac{g(m)}{m}\sum\limits_{n=1}^{+\infty}\frac{f(n)}{n^{4k-1}}\biggl(a_{k}\biggl(\frac{n}{m}\biggr) - \frac{m}{n}\biggr).
\]
Next, in view of (\ref{lab_23}), the inner sum splits into two absolutely convergent series. Thus we obtain:
\[
\frac{1}{2}\,L^{2}(2k;f) = \frac{1}{2}\sum\limits_{m=1}^{+\infty}\frac{g(m)}{m}
\biggl(\,\sum\limits_{n=1}^{+\infty}\frac{f(n)}{n^{4k-1}}\,a_{k}\biggl(\frac{n}{m}\biggr) - mL(4k;f)\biggr).
\]
Theorem is proved. $\square$
\vspace{0.5cm}

\noindent
\textsc{Corollary 1.} \emph{Under the assumptions of Theorem 4, the following identity holds:}
\[
L^{2}(2;f) =  \sum\limits_{m=1}^{+\infty}\frac{g(m)}{m}\biggl(\,2\pi\sum\limits_{n=1}^{+\infty}\frac{f(n)}{n^{3\mathstrut}\bigl(e^{2\pi n/m}-1\bigr)} - mL(4;f)+\pi L(3;f)\biggr).
\]
\textsc{Proof.} Indeed, for $k = 1$, we have
\[
L^{2}(2;f) =  \sum\limits_{m=1}^{+\infty}\frac{g(m)}{m}\biggl(\,\sum\limits_{n=1}^{+\infty}\frac{f(n)}{n^{3\mathstrut}}\,a_{1}\biggl(\frac{n}{m}\biggr) - mL(4;f)\biggr).
\]
Now it remains to note that
\[
a_{1}(w) = \pi i\ctg{(\pi iw)} = \pi \cth{(\pi w)} = \pi\biggl(1 + \frac{2}{e^{2\pi w}-1}\biggr).\quad \square
\]
In what follows, we present several particular cases of Corollary 1.
\vspace{0.5cm}

\noindent
$1^{\circ}$. Taking $g(m) = \mu(m)$, we obtain
\begin{equation*}
f(n) =
\begin{cases}
1, & \text{ if } n = 1,\\
0, & \text{ if } n > 1.
\end{cases}
\end{equation*}
Hence, $L(z;f)\equiv 1$ for any $z$. Therefore, we have
\[
\sum\limits_{m=1}^{+\infty}\frac{\mu(m)}{m}\biggl(\frac{2\pi}{e^{2\pi/m}-1} + \pi - m\biggr) = 1.
\]
If we express the left-hand side as the limit of partial sums as $M\to +\infty$ and use the fact that
\[
\lim_{M\to +\infty}\sum\limits_{1\leqslant m\leqslant M}\frac{\mu(m)}{m} = 0,
\]
we find
\[
\sum\limits_{m=1}^{+\infty}\mu(m)\biggl(\frac{2\pi/m}{e^{2\pi/m}-1}\;-\; 1\biggr) = 1.
\]

\noindent
$2^{\circ}$. Let $\nu\geqslant 1$ be an integer. Setting $g(m) = \tau_{\nu}(m)$, we obtain $f(n) = \tau_{\nu+1}(n)$ and $L(z;f) = \zeta^{\nu+1}(z)$. Hence
\[
\zeta^{2\nu+2}(2) = \sum\limits_{m=1}^{+\infty}\frac{\tau_{\nu}(m)}{m}\biggl(\,\sum\limits_{n=1}^{+\infty}\frac{\tau_{\nu+1}(n)}{n^{3\mathstrut}}\,\frac{2\pi}{e^{2\pi n/m}-1}
-m\zeta^{\nu+1}(4) + \pi\zeta^{\nu+1}(3)\biggr).
\]
In particular, for $\nu = 1$ and $\nu = 2$ we  respectively obtain
\begin{align*}
& \zeta^{4}(2) = \sum\limits_{m=1}^{+\infty}\frac{1}{m}\biggl(\,\sum\limits_{n=1}^{+\infty}\frac{\tau(n)}{n^{3\mathstrut}}\,\frac{2\pi}{e^{2\pi n/m}-1} - m\zeta^{2}(4) + \pi\zeta^{2}(3)\biggr),\\
& \zeta^{6}(2) = \sum\limits_{m=1}^{+\infty}\frac{\tau(m)}{m}\biggl(\,\sum\limits_{n=1}^{+\infty}\frac{\tau_{3}(n)}{n^{3\mathstrut}}\,\frac{2\pi}{e^{2\pi n/m}-1} - m\zeta^{3}(4) + \pi\zeta^{3}(3)\biggr).
\end{align*}

\noindent
$3^{\circ}$. Setting $g(m) = \mu^{2}(m)$, we obtain
\[
f(n) = \sum\limits_{d|n}\mu^{2}(d) = 2^{\omega(n)},\quad L(z;f) = \frac{\zeta^{2}(z)}{\zeta(2z)}\quad \text{for}\quad \RRe z > 1,
\]
and hence,
\[
\frac{\zeta^{2}(2)}{\zeta(4)} = \sum\limits_{m=1}^{+\infty}\frac{\mu^{2}(m)}{m}\biggl(\,\sum\limits_{n=1}^{+\infty}\frac{2^{\omega(n)}}{n^{3\mathstrut}}\,\frac{2\pi}{e^{2\pi n/m}-1} -
m\,\frac{\zeta^{2}(4)}{\zeta(8)} + \pi\,\frac{\zeta^{2}(3)}{\zeta(6)}\biggr).
\]
\noindent
$4^{\circ}$. Setting
\begin{equation*}
g(m) =
\begin{cases}
\mu(d), & \text{ if } m = d^{2},\\
0, & \text{ otherwise},
\end{cases}
\quad \text{we get}\quad
f(n) =
\begin{cases}
1, & \text{ if } n \text{ is squarefree},\\
0, & \text{ otherwise},
\end{cases}
\end{equation*}
Thus
\[
f(n) = \mu^{2}(n),\quad L(z;f) = \frac{\zeta(z)}{\zeta(2z)},\quad \text{for}\quad \RRe z > 1,
\]
and hence,
\[
\biggl(\frac{\zeta(2)}{\zeta(4)}\biggr)^{2} = \sum\limits_{m=1}^{+\infty}\frac{\mu(m)}{m^{2\mathstrut}}\biggl(\,\sum\limits_{n=1}^{+\infty}\frac{\mu^{2}(n)}{n^{3\mathstrut}}\,\frac{2\pi}{e^{2\pi n/m^{2}}-1}
-m^{2}\,\frac{\zeta(4)}{\zeta(8)} + \pi\,\frac{\zeta(3)}{\zeta(6)}\biggr).
\]

\noindent
$5^{\circ}$. Setting $g(m) = \tau(m^{2})$ and using the identities
\[
\frac{\zeta^{3}(z)}{\zeta(2z)} = \sum\limits_{m=1}^{+\infty}\frac{\tau(m^{2})}{m^{z}},\quad
\frac{\zeta^{4}(z)}{\zeta(2z)} = \sum\limits_{m=1}^{+\infty}\frac{\tau^{2}(m)}{m^{z}},\quad \RRe z > 1,
\]
we conclude that
\[
f(n) = \tau^{2}(n),\quad L(z;f) = \frac{\zeta^{4}(z)}{\zeta(2z)}\quad\text{for}\quad \RRe z > 1
\]
and hence,
\[
\frac{\zeta^{8}(2)}{\zeta^{2}(4)} = \sum\limits_{m=1}^{+\infty}\frac{\tau(m^{2})}{m}
\biggl(\,\sum\limits_{n=1}^{+\infty}\frac{\tau^{2}(n)}{n^{3\mathstrut}}\,\frac{2\pi}{e^{2\pi n/m}-1}
-m\,\frac{\zeta^{4}(4)}{\zeta(8)} + \pi\,\frac{\zeta^{4}(3)}{\zeta(6)}\biggr).
\]

\noindent
$6^{\circ}$. Setting $g(m) = \lambda(m)$, we obtain
\begin{equation*}
f(n) =
\begin{cases}
1, & \text{if } n \text{ is a perfect square},\\
0, & \text{otherwise},
\end{cases}
\quad L(z;f) = \zeta(2z) \quad\text{for}\quad \RRe z > \frac{1}{2},
\end{equation*}
and hence,
\[
\zeta^{2}(4) = \sum\limits_{m=1}^{+\infty}\frac{\lambda(m)}{m}
\biggl(\,\sum\limits_{n=1}^{+\infty}\frac{1}{n^{6\mathstrut}}\,\frac{2\pi}{e^{2\pi n^{2\mathstrut}/m}-1}
-m\,\zeta(8) + \pi\,\zeta(6)\biggr).
\]

\noindent
$7^{\circ}$. Setting $g(m) = \mu(m)/m$, we obtain
\[
f(n) = \frac{\varphi(n)}{n},\quad L(z;f) = \frac{\zeta(z)}{\zeta(z+1)}\quad\text{for}\quad \RRe z > 1
\]
and thus,
\[
\biggl(\frac{\zeta(2)}{\zeta(3)}\biggr)^{2} = \sum\limits_{m=1}^{+\infty}\frac{\mu(m)}{m^{2\mathstrut}}\biggl(\,\sum\limits_{n=1}^{+\infty}\frac{\varphi(n)}{n^{4\mathstrut}}\,\frac{2\pi}{e^{2\pi n/m}-1} -m\,\frac{\zeta(4)}{\zeta(5)} + \pi\,\frac{\zeta(3)}{\zeta(4)}\biggr).
\]

\noindent
$8^{\circ}$. Setting $g(m)= \Lambda(m)$, we obtain $f(n) = \ln{n}$, $L(z;f) = -\zeta'(z)$; hence we have
\[
\bigl(\zeta'(2)\bigr)^{2} = \sum\limits_{m=1}^{+\infty}\frac{\Lambda(m)}{m}
\biggl(\,\sum\limits_{n=1}^{+\infty}\frac{\ln{n}}{n^{3\mathstrut}}\,\frac{2\pi}{e^{2\pi n/m}-1}
+m\,\zeta'(4) - \pi\,\zeta'(3)\biggr).
\]

\noindent
$9^{\circ}$. Suppose that $k\geqslant 1$, and let $\Lambda_{k}$ be the generalized von Mangoldt function defined by the following relation:
\[
\sum\limits_{m=1}^{+\infty}\frac{\Lambda_{k}(m)}{m^{z}} = (-1)^{k}\,\frac{\zeta^{(k)}(z)}{\zeta(z)}.
\]
Taking $g(m) = \Lambda_{k}(m)$, we obtain $f(n) = (-1)^{k}(\ln{n})^{k}$, $L(z;f) = (-1)^{k}\zeta^{(k)}(z)$; therefore
\[
\bigl(\zeta^{(k)}(2)\bigr)^{2} = (-1)^{k}\sum\limits_{m=1}^{+\infty}\frac{\Lambda_{k}(m)}{m}
\biggl(\,\sum\limits_{n=1}^{+\infty}\frac{(\ln{n})^{k}}{n^{3\mathstrut}}\,\frac{2\pi}{e^{2\pi n/m}-1}
-m\,\zeta^{(k)}(4) + \pi\,\zeta^{(k)}(3)\biggr).
\]

\noindent
$10^{\circ}$. Let $\chi_{4}$ be the non-principal Dirichlet character modulo $4$, given by
\begin{equation*}
\chi_{4}(m) =
\begin{cases}
1, & \text{if } m\equiv 1 \pmod{4},\\
-1, & \text{if } m\equiv 3 \pmod{4},\\
0, & \text{if } m\equiv 0 \pmod{2}.
\end{cases}
\end{equation*}
Then, setting $g(m) = \chi_{4}(m)$, we get $f(n) = r(n)/4$, where $r(n)$ denotes the number of representations of $n$ as a sum of two squares of integers. Then
\[
L(z;f) = \frac{1}{4}\,\zeta(z)L(z),\quad\text{where}\quad L(z) = L(z;\chi_{4}) = \sum\limits_{m=1}^{+\infty}\frac{\chi_{4}(m)}{m^{z}}
\]
is the Dirichlet $L$-function associated with $\chi_{4}$. After straightforward calculations, we obtain
\[
\zeta^{2}(2)\beta^{2} = 4\sum\limits_{m=1}^{+\infty}\frac{\chi_{4}(m)}{m}
\biggl(\,\sum\limits_{n=1}^{+\infty}\frac{r(n)}{n^{3\mathstrut}}\,\frac{2\pi}{e^{2\pi n/m}-1}
-m\,\zeta(4)L(4) + \pi\,\zeta(3)L(3)\biggr),
\]
where
\[
\beta = L(2) = \sum\limits_{n=0}^{+\infty}\frac{(-1)^{n}}{(2n+1)^{2}}\quad \text{is Catalan constant.}
\]

\noindent
$11^{\circ}$. Let $g(m) = c_{m}(a)$, where
\[
c_{m}(a) = \sum\limits_{\substack{\nu = 1 \\ (\nu,m)=1}}^{m}e^{2\pi i\,\frac{\scriptstyle \nu a}{\scriptstyle m}}\quad \text{is Ramanujan's sum}
\]
(see \cite{Ramanujan_1918}; here and below, we suppose that $a\geqslant 1$ is a fixed integer). Then
\begin{equation*}
f(n) = \sum\limits_{m|n}g(m) = \sum\limits_{m|n}c_{m}(a) =
\begin{cases}
n, & \text{if } n|a,\\
0, & \text{otherwise}.
\end{cases}
\end{equation*}
Hence,
\[
L(z;f) = \sum\limits_{n|a}\frac{n}{n^{z}} = \sigma_{1-z}(a) = \frac{\sigma_{z-1}(a)}{a^{z-1}}
\]
and, in particular,
\[
L(2;f) = \frac{\sigma(a)}{a},\quad L(3;f) = \frac{\sigma_{2}(a)}{a^{2\mathstrut}},\quad L(4;f) = \frac{\sigma_{3}(a)}{a^{3\mathstrut}}.
\]
Thus we find
\[
\biggl(\frac{\sigma(a)}{a}\biggr)^{2} = \sum\limits_{m=1}^{+\infty}\frac{c_{m}(a)}{m}\biggl(\,\sum\limits_{d|a}\frac{2\pi}{d^{2}(e^{2\pi d/m}-1)} - m\,\frac{\sigma_{3}(a)}{a^{3\mathstrut}}+\pi\,\frac{\sigma_{2}(a)}{a^{2\mathstrut}}\biggr).
\]
Note that
\[
\sum\limits_{m=1}^{+\infty}\frac{c_{m}(a)}{m} = 0\quad\text{for any } a\geqslant 1.
\]
Therefore,
\[
\sigma^{2}(a) = a^{2}\sum\limits_{m=1}^{+\infty}c_{m}(a)\biggl(\,\sum\limits_{d|a}\frac{1}{d^{2\mathstrut}}\,\frac{2\pi/m}{e^{2\pi d/m}-1}-\frac{\sigma_{3}(a)}{a^{3\mathstrut}}\biggr).
\]
In the case $a=1$, using the fact that $c_{m}(1) = \mu(m)$, we obtain Cauchy-Lerch-Ramanujan identity from Case $1^{\circ}$.
\vspace{0.5cm}

\noindent
\textsc{Theorem 5.} \emph{Under the assumptions of Theorem 4, for any $k\geqslant 2$, the following identity holds}
\[
L^{2}(2k-1;f) = \sum\limits_{m=1}^{+\infty}\frac{g(m)}{m}\sum\limits_{n=1}^{+\infty}\frac{f(n)}{n^{4k-3}}\,b_{k}\biggl(\frac{n}{m}\biggr).
\]
\textsc{Proof.} Setting $\ell = 1$ in Corollary 2 of Lemma 3, we obtain
\[
\frac{qn^{2k-2}}{q^{2k\mathstrut}+n^{2k\mathstrut}} = \frac{1}{2k}\sum\limits_{r=0}^{2k-1}\varepsilon_{r}^{2}\biggl(\frac{1}{q} - \frac{1}{q+n\varepsilon_{r}}\biggr)
\]
and therefore
\[
\sum\limits_{q=1}^{+\infty}\frac{n^{2k-2}\,qf(q)}{q^{2k\mathstrut}+n^{2k\mathstrut}} = \frac{1}{2k}\sum\limits_{q=1}^{+\infty}f(q)
\sum\limits_{r=0}^{2k-1}\varepsilon_{r}^{2}\biggl(\frac{1}{q} - \frac{1}{q+n\varepsilon_{r}}\biggr).
\]
Using the inequality
\[
|q+n\varepsilon_{r}|\geqslant (q+n)\sin{\frac{\pi}{4k}}
\]
and following word for word the proof of Theorem 4, we obtain
\begin{multline*}
\sum\limits_{q=1}^{+\infty}\frac{n^{2k-2}\,qf(q)}{q^{2k\mathstrut}+n^{2k\mathstrut}} = \frac{1}{2k}\sum\limits_{r=0}^{2k-1}\varepsilon_{r}^{2}\sum\limits_{m,\ell=1}^{+\infty}g(m)\biggl(\frac{1}{m\ell} - \frac{1}{m\ell + n\varepsilon_{r}}\biggr) = \\
= \frac{1}{2k}\sum\limits_{r=0}^{2k-1}\varepsilon_{r}^{2}\sum\limits_{m=1}^{+\infty}\frac{g(m)}{m}\sum\limits_{\ell =1}^{+\infty}\biggl(\frac{1}{\ell} - \frac{1}{\ell + n\varepsilon_{r}/m}\biggr) = \frac{1}{2k}\sum\limits_{r=0}^{2k-1}\varepsilon_{r}^{2}\sum\limits_{m=1}^{+\infty}\frac{g(m)}{m}\biggl(\psi\biggl(\frac{n\varepsilon_{r}}{m}\biggr)+\gamma + \frac{m}{n\varepsilon_{r}}\biggr) =\\
= \frac{1}{2}\sum\limits_{m=1}^{+\infty}\frac{g(m)}{m}\,\frac{1}{k}\sum\limits_{r=0}^{2k-1}\varepsilon_{r}^{2}\biggl(\psi\biggl(\frac{n\varepsilon_{r}}{m}\biggr)+\gamma + \frac{m}{n\varepsilon_{r}}\biggr) = \frac{1}{2}\sum\limits_{m=1}^{+\infty}\frac{g(m)}{m}\,b_{k}\biggl(\frac{n}{m}\biggr).
\end{multline*}
Multiplying both sides by $f(n)\cdot n^{-(4k-3)}$ and summing over $n$, we find
\begin{multline*}
\frac{1}{2}\sum\limits_{n=1}^{+\infty}\frac{f(n)}{n^{4k-3}}\sum\limits_{m=1}^{+\infty}\frac{g(m)}{m}\,b_{k}\biggl(\frac{n}{m}\biggr) =
\frac{1}{2}\sum\limits_{m=1}^{+\infty}\frac{g(m)}{m}\sum\limits_{n=1}^{+\infty}\frac{f(n)}{n^{4k-3}}\,b_{k}\biggl(\frac{n}{m}\biggr) = \\
\frac{1}{2}\sum\limits_{n=1}^{+\infty}\frac{f(n)}{n^{4k-3}}\sum\limits_{q=1}^{+\infty}\frac{n^{2k-2}\,qf(q)}{q^{2k\mathstrut}+n^{2k\mathstrut}} = \sum\limits_{q,n=1}^{+\infty}\frac{f(q)f(n)}{(qn)^{2k-1}}\,\frac{q^{2k}}{q^{2k\mathstrut} + n^{2k\mathstrut}} =\\
= \frac{1}{2}\biggl(\,\sum\limits_{n=1}^{+\infty}\frac{f(n)}{n^{2k-1}}\biggr)^{2} = \frac{1}{2}\,L^{2}(2k-1;f).\quad \square
\end{multline*}
\textsc{Corollary.} \emph{Let $\varepsilon = e^{\pi i/4}$, $\omega = e^{\pi i/6}$. Then, under the assumptions of Theorem 4, the following identities hold:}
\begin{align*}
L^{2}(3;f) = & \frac{i}{2}\sum\limits_{m=1}^{+\infty}\frac{g(m)}{m}\sum\limits_{n=1}^{+\infty}\frac{f(n)}{n^{5}}\,\biggl(\psi\Bigl(\frac{n\varepsilon}{m}\Bigr) + \psi\Bigl(-\frac{n\varepsilon}{m}\Bigr) - \psi\Bigl(\frac{n\overline{\varepsilon}}{m}\Bigr) - \psi\Bigl(-\frac{n\overline{\varepsilon}}{m}\Bigr)\biggr),\\
L^{2}(5;f) = & \frac{1}{3}\sum\limits_{m=1}^{+\infty}\frac{g(m)}{m}\sum\limits_{n=1}^{+\infty}\frac{f(n)}{n^{9}}\,\biggl\{\omega^{2}\biggl(\psi\Bigl(\frac{n\omega}{m}\Bigr) + \psi\Bigl(-\frac{n\omega}{m}\Bigr)\biggr) +\overline{\omega}^{2}\biggl(\psi\Bigl(\frac{n\overline{\omega}}{m}\Bigr) + \psi\Bigl(-\frac{n\overline{\omega}}{m}\Bigr)\biggr) -\\
- & \psi\Bigl(\frac{ni}{m}\Bigr) - \psi\Bigl(-\frac{ni}{m}\Bigr)\biggr\}.
\end{align*}

\textsc{Theorem 6.} \emph{Under the assumptions of Theorem 4, for any $k\geqslant 1$, the following identity holds:}
\[
\frac{1}{2}\,L^{2}(2k+1;f) = \sum\limits_{m=1}^{+\infty}\frac{g(m)}{m}\biggl(\,\sum\limits_{n=1}^{+\infty}\frac{f(n)}{n^{4k+1}}\,c_{k}\Bigl(\frac{n}{m}\Bigr)- L(4k+2;f)\biggr).
\]
\textsc{Proof} is based on the expansion
\[
\frac{n^{2k}}{q^{2k+1\mathstrut} + n^{2k+1\mathstrut}} = \frac{1}{2k+1}\sum\limits_{r=0}^{2k}\omega_{r}\biggl(\frac{1}{q} - \frac{1}{q-n\omega_{r}}\biggr)
\]
(see Corollary 1 to Lemma 4) and follows word for word the proofs of Theorems 4, 5.
\vspace{0.5cm}

\textsc{Corollary.} \emph{Let $\varpi = e^{\pi i/3}$. Then, under the assumptions of Theorem 4, the following identity holds:}
\[
\frac{1}{2}\,L^{2}(3;f) = \sum\limits_{m=1}^{+\infty}\frac{g(m)}{m}\biggl(\frac{1}{3}\sum\limits_{n=1}^{+\infty}\frac{1}{n^{5}}\biggl\{\varpi\psi\Bigl(-\frac{n\varpi}{m}\Bigr) +
\overline{\varpi}\psi\Bigl(-\frac{n\overline{\varpi}}{m}\Bigr) - \psi\Bigl(\frac{n}{m}\Bigr)\biggr\} - L(6;f)\biggr).
\]
\vspace{0.5cm}

\textbf{\S 5. The limiting behavior of the coefficients $\boldsymbol{a_{k}(w)}$, $\boldsymbol{b_{k,\ell}(w)}$, and $\boldsymbol{c_{k}(w)}$.}
\vspace{0.3cm}

The goal of this section is to study the behavior of the coefficients $a_{k}(w)$, $b_{k,\ell}(w)$, and $c_{k}(w)$ for fixed $w$ and $\ell$ as $k$ tends to infinity.
\vspace{0.3cm}

\textsc{Theorem 7.} \emph{Suppose that $w>0$ and an integer $\ell\geqslant 1$ are fixed. Then, as $k\to +\infty$, the following relations hold:}
\[
a_{k}(w)\to \mathfrak{a}(w),\quad b_{k,\ell}(w)\to \mathfrak{b}_{\ell}(w),\quad c_{k}(w)\to \mathfrak{c}(w),
\]
\emph{where}
\begin{equation*}
\mathfrak{a}(w) =
\begin{cases}
2, & \textit{for integer }\; w,\\[6pt]
\displaystyle \frac{2[w]+1}{w}, & \textit{otherwise};
\end{cases}
\quad \mathfrak{c}(w) =
\begin{cases}
\displaystyle \frac{w+1/2}{w}, & \textit{for integer }\; w,\\[6pt]
\displaystyle \frac{[w]+1}{w}, & \textit{otherwise}.
\end{cases}
\end{equation*}
\begin{equation*}
\mathfrak{b}_{\ell}(w) = \frac{2}{w^{\ell+1}}\biggl(\frac{B_{\ell+1}([w])-B_{\ell+1}}{\ell+1} + \kappa[w]^{\ell}\biggr),\quad
\kappa =
\begin{cases}
1/2, & \textit{for integer }\; w,\\[6pt]
1, & \textit{otherwise},
\end{cases}
\end{equation*}
\emph{and $B_{j}(x)$ is $j$th Bernoulli polynomial.}
\vspace{0.5cm}

\textsc{Proof.} To establish the first relation, we set
\[
g(u) = \pi e^{2\pi iu}\ctg{\bigl(\pi w e^{2\pi iu}\bigr)},\quad f(x) = g\biggl(\frac{x+1/2}{2k}\biggr)
\]
and note that
\[
a_{k}(w) = \frac{1}{k}\sum\limits_{r=0}^{k-1}\pi\varepsilon_{r}\ctg{\bigl(\pi w\varepsilon_{r}\bigr)} = \frac{1}{k}\sum\limits_{r=0}^{k-1}f(r).
\]
Assuming that $0<\varepsilon< 1/2$, we obtain
\[
a_{k}(w) = \frac{1}{k}\sum\limits_{\xi< r\leqslant \eta} f(r),\quad\text{where}\quad \xi = -\frac{1}{2}+\varepsilon,\quad \eta = k-\frac{1}{2}-\varepsilon.
\]
By the Euler-Maclaurin summation formula,
\begin{multline}\label{lab_25}
a_{k}(w) = \\
= \frac{1}{k}\biggl(\,\int_{\xi}^{\eta}f(x)dx + \int_{\xi}^{\eta}\sigma(x)f''(x)dx + \varrho(\eta)f(\eta) - \varrho(\xi)f(\xi) + \sigma(\xi)f'(\xi) - \sigma(\eta)f'(\eta)\biggr).
\end{multline}
Since $g^{(r)}(u+1/2) = g^{(r)}(u)$ and $g^{(r)}(-u) = (-1)^{r}\bar{g}^{(r)}(u)$ for any $r\geqslant 0$, then, setting $\delta = \varepsilon/(2k)$ we get
\begin{multline*}
f(\xi) = g(\delta),\quad f(\eta) = g(1/2-\delta) = g(-\delta) = \bar{g}(\delta),\\
f'(\xi) = \frac{1}{2k}\,g'(\delta),\quad f'(\eta) = \frac{1}{2k}\,g'(1/2-\delta) = -\frac{1}{2k}\,\bar{g}'(\delta).
\end{multline*}
Further,
\[
\varrho(\xi) = -\varepsilon,\quad \varrho(\eta) = \varepsilon,\quad \sigma(\xi) = \sigma(\eta) = \frac{1}{2}\biggl(\frac{1}{4}-\varepsilon^{2}\biggr).
\]
Hence, the boundary terms in (\ref{lab_25}) equal
\begin{equation}\label{lab_26}
4\delta \RRe g(\delta) + \frac{1}{2k^{2}}\biggl(\frac{1}{4}-\varepsilon^{2}\biggr)\RRe g'(\delta).
\end{equation}
Next, the first integral in (\ref{lab_25}) coincides with
\begin{equation}\label{lab_27}
\frac{1}{k}\int_{\xi}^{\eta}g\biggl(\frac{x+1/2}{2k}\biggr)dx = 2\int_{\delta}^{1/2-\delta}g(u)du,
\end{equation}
and the second one is transformed as follows:
\begin{multline}\label{lab_28}
\frac{1}{k}\cdot\frac{1}{(2k)^{2}}\int_{\xi}^{\eta}\sigma(x)g''\biggl(\frac{x+1/2}{2k}\biggr)dx = \frac{1}{2k^{2}}\int_{\delta}^{1/2-\delta}\sigma\bigl(2ku-1/2\bigr)g''(u)du = \\
= \frac{1}{2k^{2}}\biggl(\,\int_{\delta}^{1/4} + \int_{1/4}^{1/2-\delta}\biggr)\sigma\bigl(2ku-1/2\bigr)g''(u)du = \\
= \frac{1}{2k^{2}}\int_{\delta}^{1/4}\sigma\bigl(2ku-1/2\bigr)\bigl(g''(u) + g''(1/2-u)\bigr)du = \\
= \frac{1}{k^{2}}\int_{\delta}^{1/4}\sigma\bigl(2ku-1/2\bigr)\RRe{g''(u)}\,du.
\end{multline}
Straightforward calculations yield:
\begin{align*}
g'(u) & = 2\pi^{2}i e^{2\pi i u}\biggl(\ctg{\bigl(\pi w e^{2\pi iu}\bigr)} - \frac{\pi w e^{2\pi iu}}{\sin^{2}(\pi w e^{2\pi iu})}\biggr),\\[6pt]
g''(u) & = -4\pi^{3}e^{2\pi iu}\biggl(\ctg{\bigl(\pi w e^{2\pi iu}\bigr)} - \frac{3\pi w e^{2\pi iu}}{\sin^{2}(\pi w e^{2\pi iu})} + \frac{2(\pi w e^{2\pi iu})^{2}\cos{(\pi w e^{2\pi iu})}}{\sin^{3}(\pi w e^{2\pi iu})}\biggr).\\
\end{align*}
\begin{figure}[h]
\begin{minipage}[h]{1.0\linewidth}
\center{\includegraphics[width=1.0\linewidth]{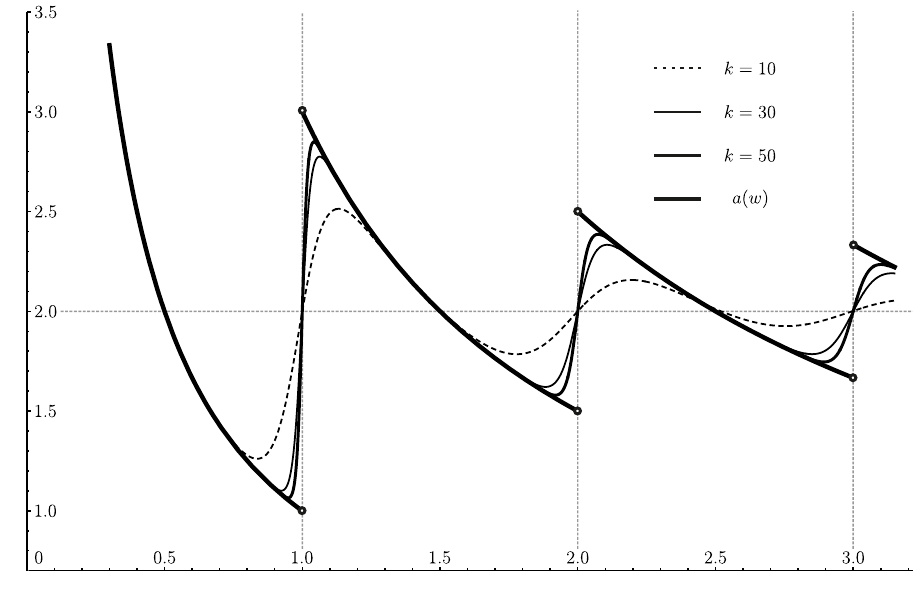}}\\
\end{minipage}
\begin{center}
\emph{Plots of the functions $\mathfrak{a}(w)$ and $a_{k}(w)$ for $k = 10, 30, 50$.}
\end{center}
\end{figure}
For non-integer $w$, the functions $g^{(r)}(u)$, $r = 0,1,2$, have no singularities at $u=0$ and $u = 1/2$. In addition, the value
\[
g'(0) = 2\pi^{2}i\biggl(\ctg{(\pi w)} - \frac{\pi w}{\sin^{2}(\pi w)}\biggr)
\]
is purely imaginary. Hence, as $\varepsilon\to 0$, the sum (\ref{lab_26}) tends to $(8k^{2})^{-1}\RRe{g'(0)} = 0$, while the right-hand sides of (\ref{lab_27}) and (\ref{lab_28}) tend to the integrals
\[
\mathfrak{a}(w) = 2\int_{0}^{1/2}g(u)du = \int_{-1/2}^{1/2}g(u)du,\quad \alpha_{k}(w) = \int_{0}^{1/4}\sigma(2ku-1/2)\RRe{g''(u)}du,
\]
respectively. Thus we obtain
\begin{equation}\label{lab_29}
a_{k}(w)= \mathfrak{a}(w) + \frac{\alpha_{k}(w)}{k^{2}}.
\end{equation}
Suppose now that $w = m\geqslant 1$ is an integer. Then, as $\xi\to 0$, the following expansions hold:
\begin{align*}
& g(\xi) = \frac{1}{2\pi im}\,\frac{1}{\xi} + \frac{1}{2m} - \frac{\pi i}{6m}\,\bigl((2\pi m)^{2}-1\bigr)\xi + 2\pi^{4}m\xi^{2} + O(|\xi|^{3}),\\
& g'(\xi) = -\,\frac{1}{2\pi i m}\,\frac{1}{\xi^{2}} - \frac{\pi i}{6m}\,\bigl((2\pi m)^{2}-1\bigr) + 4\pi^{4}m\xi + O(|\xi|^{2}),\\
& g''(\xi) = \frac{1}{\pi im}\,\frac{1}{\xi^{3}} + 4\pi^{4}m + O(|\xi|).
\end{align*}
Therefore,
\[
4\delta\RRe{g(\delta)} = \frac{2\delta}{m} + O(\delta^{3}),\quad \RRe{g'(\delta)} = 4\pi^{4}m\delta + O(\delta^{2}),\quad \RRe{g''(\delta)} = 4\pi^{4}m + O(\delta).
\]
Hence, the sum (\ref{lab_26}) tends to zero as $\delta\to 0$, while the integrals (\ref{lab_27}) and (\ref{lab_28}) remain bounded. Thus, the relation (\ref{lab_29}) still holds for
\[
\mathfrak{a}(w) = 2\,\text{v.p.}\,\int_{0}^{1/2}g(u)\,du.
\]
For fixed $w$, in both cases we have
\[
|\alpha_{k}(w)|\leqslant \frac{1}{8}\int_{0}^{1/4}\bigl|\RRe{g''(u)}\bigr|\,du,
\]
where the right-hand side does not depend on $k$. Hence, $a_{k}(w)\to \mathfrak{a}(w)$ as $k\to +\infty$.

Now, let us calculate $\mathfrak{a}(w)$. If $w$ is not an integer, then, setting $z = we^{2\pi iu}$, we obtain
\[
\mathfrak{a}(w) = \frac{\mathfrak{A}(w)}{w},\quad \mathfrak{A}(w) = \frac{1}{2\pi i}\int_{\mathcal{C}_{w}}\pi\ctg{(\pi z)}\,dz,
\]
where $\mathcal{C}_{w}$ denotes the circle of radius $w$ centered  at the origin. If $[w] = m$, then the integral $\mathfrak{A}(w)$ is equal to the sum of residues of the function $\pi\ctg{(\pi z)}$ at the points $z = 0,\pm 1,\pm 2,\ldots, \pm m$. Since each of these residues equals $1$, we have
\[
\mathfrak{A}(w) = 2m+1 = 2[w]+1,\quad \mathfrak{a}(w) = \frac{2[w]+1}{w}.
\]
If $w = m$ is integer then
\[
\mathfrak{a}(w) = \frac{\mathfrak{A}(w)}{w},\quad \mathfrak{A}(w) = \frac{\text{v.p.}}{2\pi i}\int_{\mathcal{C}_{w}}\pi\ctg{(\pi z)}\,dz.
\]
The symbol ``v.p.'' denotes the limit of the integrals over $\mathcal{C}_{w}$ with excluded arcs of length $\varepsilon$ centered at the points $\pm m$, as $\varepsilon\to 0$.
An straightforward calculation yields $\mathfrak{A}(w) = 2m = 2w$ and $\mathfrak{a}(w) = 2$.

The proof of the limit relation for $b_{k,\ell}(w)$ follows along the same lines. Setting
\[
g(u) = \bigl(e^{2\pi iu}\bigr)^{\ell+1}\bigl\{\psi(we^{2\pi iu}) + (-1)^{\ell+1}\psi(-we^{2\pi iu})\bigr\},\quad f(x) = g\biggl(\frac{x+1/2}{2k}\biggr),
\]
we obtain $g^{(r)}(u+1/2) = g^{(r)}(u)$, $g^{(r)}(-u) = (-1)^{r}\,\overline{g}^{(r)}(u)$.
Keeping the notation for $\xi$, $\eta$, $\delta$ and using the Euler-Maclaurin summation formula, we obtain
\begin{multline}\label{lab_30}
(-1)^{\ell+1}b_{k,\ell}(w) = \frac{1}{k}\sum\limits_{r=0}^{2k-1}\varepsilon_{r}^{\ell+1}\psi(\varepsilon w) =
\frac{1}{k}\sum\limits_{r=0}^{k-1}\varepsilon_{r}^{\ell+1}\bigl\{\psi(\varepsilon_{r}w) + (-1)^{\ell+1}\psi(-\varepsilon_{r}w)\bigr\} = \\
= \frac{1}{k}\sum\limits_{\xi<r\leqslant \eta}f(r) = 2\int_{\delta}^{1/2-\delta}g(u)du + \frac{1}{k^{2}}\int_{\delta}^{1/4}\sigma\bigl(2ku-1/2\bigr)\RRe{g''(u)}du + \\
+ 2\delta\RRe{g(\delta)} + \frac{1}{2k^{2}}\biggl(\frac{1}{4}-\varepsilon^{2}\biggr)\RRe{g'(\delta)},
\end{multline}
If $w$ is not an integer, then the functions $g^{(r)}(u)$, $r=0,1,2$ have no singularities at $u=0$ and $u = 1/2$, and the value
\[
g'(0) = 2\pi i\bigl((\ell+1)\bigl\{\psi(w)+(-1)^{\ell+1}\psi(-w)\bigr\}+w\bigl\{\psi'(w)-(-1)^{\ell+1}\psi'(-w)\bigr\}\bigr)
\]
is purely imaginary. Letting $\varepsilon\to 0$, we obtain
\begin{multline}\label{lab_31}
b_{k,\ell}(w) = \mathfrak{b}_{\ell}(w) + \frac{\beta_{k,\ell}(w)}{k^{2}},\\
\mathfrak{b}_{\ell}(w) = 2(-1)^{\ell+1}\int_{0}^{1/2}g(u)\,du = (-1)^{\ell+1}\int_{-1/2}^{1/2}g(u)\,du,\\
 \beta_{k,\ell}(w) = (-1)^{\ell+1}\int_{0}^{1/4}\sigma(2ku-1/2)\RRe{g''(u)}\,du.
\end{multline}
For an integer $w = m\geqslant 1$, we have
\begin{align*}
& g(\xi) = \frac{1}{2\pi im}\frac{1}{\xi} + A_{m} + iB_{m}\xi + C_{m}\xi^{2} + O(|\xi|^{3}),\\
& g'(\xi) = -\frac{1}{2\pi im}\frac{1}{\xi^{2}} + iB_{m} + 2C_{m}\xi + O(|\xi|^{2}),\\
& g''(\xi) = \frac{1}{\pi im}\frac{1}{\xi^{3}} + 2C_{m} + O(|\xi|),
\end{align*}
where $A_{m}, B_{m}, C_{m}$ are real coefficients depending only on $m$ and $\ell$. Hence,
\[
2\delta\RRe{g(\delta)} = 2\delta A_{m} + O(\delta^{3}),\quad \RRe{g'(\delta)} = 2C_{m}\delta + O(\delta^{2}),\quad \RRe{g''(\delta)} = 2C_{m} + O(\delta),
\]
as $\delta \to 0$. Therefore, the boundary terms in (\ref{lab_30}) vanish as $\varepsilon\to 0$ and the formula (\ref{lab_31}) remains valid with
\[
\mathfrak{b}_{\ell}(w) = 2(-1)^{\ell+1}\,\text{v.p.}\int_{0}^{1/2}g(u)du.
\]
Taking the limit $k\to +\infty$, we obtain $b_{k,\ell}(w)\to \mathfrak{b}_{\ell}(w)$ in both cases.

For non-integral $w$, we have
\begin{equation}\label{lab_32}
\mathfrak{b}_{\ell}(w) = \frac{\mathfrak{B}_{\ell}(w)}{w^{\ell+1}},\quad \mathfrak{B}_{\ell}(w) = \frac{(-1)^{\ell+1}}{2\pi i}\int_{\mathcal{C}_{w}}z^{\ell}\bigl(\psi(z)+(-1)^{\ell+1}\psi(-z)\bigr)dz.
\end{equation}
Suppose that $r\geqslant 0$ is an integer. Then the residue of the function $z^{\ell}\psi(z)$ at the point $z = -r$ equals $(-1)^{\ell+1}r^{\ell}$ and the residue of the function $z^{\ell}\psi(-z)$ at the point $z = r$ is $r^{\ell}$. Hence,
\[
\mathfrak{B}_{\ell}(w) = (-1)^{\ell+1}\biggl(\sum\limits_{r=1}^{m}(-1)^{\ell+1}r^{\ell}+(-1)^{\ell+1}\sum\limits_{r=1}^{m}r^{\ell}\biggr) = 2\sum\limits_{r=1}^{m}r^{\ell}.
\]
If $w = m\geqslant 1$ is an integer, then (\ref{lab_32}) holds true
\[
\mathfrak{B}_{\ell}(w) =  (-1)^{\ell+1}\cdot\frac{\text{v.p.}}{2\pi i}\int_{\mathcal{C}_{w}}z^{\ell}\bigl(\psi(z)+(-1)^{\ell+1}\psi(-z)\bigr)dz = \sum\limits_{r=1}^{m}r^{\ell} + \sum\limits_{r=1}^{m-1}r^{\ell} = 2\sum\limits_{r=1}^{m}r^{\ell} - m^{\ell}.
\]
It remains to observe that
\[
\sum\limits_{r=1}^{m}r^{\ell} = \frac{1}{\ell+1}\sum\limits_{j=0}^{\ell}\binom{\ell+1}{j}B_{j}m^{\ell+1-j} + m^{\ell} = \frac{B_{\ell+1}(m)-B_{\ell+1}}{\ell+1} + m^{\ell}.
\]
To prove the final assertion, we set
\begin{multline*}
g(u) = e^{2\pi iu}\psi\bigl(-we^{2\pi iu}\bigr),\quad f(x) = g\biggl(\frac{x+1/2}{2k+1}\biggr),\quad \delta = \frac{\varepsilon}{2k+1},\\
\xi = -\frac{1}{2}+\varepsilon,\quad \eta = 2k+\frac{1}{2}-\varepsilon.
\end{multline*}
Then
\[
c_{k}(w) = \frac{1}{2k+1}\sum\limits_{\xi<r\leqslant \eta}\omega_{r}\psi(-\omega_{r}w) = \frac{1}{2k+1}\sum\limits_{\xi<r\leqslant \eta}f(r).
\]
Using the Euler-Maclaurin summation formula and taking the limit as $\varepsilon\to 0$, we obtain the relation
\[
c_{k}(w) = \mathfrak{c}(w) + \frac{\gamma_{k}(w)}{(2k+1)^{2}}
\]
where
\[
\mathfrak{c}(w) = \text{v.p.}\int_{0}^{1}g(u)du,\quad \gamma_{k}(w) = 2\int_{0}^{1/2}\sigma((2k+1)u-1/2)\RRe{g''(u)}du
\]
and
\[
\mathfrak{c}(w) = \frac{\mathfrak{C}(w)}{w},\quad \mathfrak{C}(w) = \frac{\text{v.p.}}{2\pi i}\int_{C_{w}}\psi(-z)dz
\]
(for non-integer $w$, the symbol ``v.p.'' can be omitted). A simple calculation shows that
$\mathfrak{C}(w) = [w]+1$ for non-integral $w$ and $\mathfrak{C}(w) = m+1/2$ for an integer $w = m$. Thus, the assertion follows.
\vspace{0.5cm}

\renewcommand{\refname}{\normalsize{References}}

\end{document}